\documentclass[12pt]{article}
\usepackage{latexsym}
\usepackage{amsfonts}
\usepackage{amsmath}
\usepackage{amssymb}
\usepackage{url}

\makeatletter
\@addtoreset{equation}{section}

\ifx\reset@font\undefined
                \let\reset@font=\relax
\fi
\def\section{\@startsection {section}{1}{\z@}{-3.5ex plus-1ex minus
                  -.2ex}{2.3ex plus.2ex}{\reset@font\large\bf}}
\def\subsection{\@startsection{subsection}{2}{\z@}{-3.25ex plus-1ex
                  minus-.2ex}{-.1em}{\reset@font\large\bf}}
\def\subsubsection{\@startsection{subsubsection}{3}{\z@}{-3.25ex plus
               -1ex minus-.2ex}{-.1em}{\reset@font\normalsize\bf}}
\makeatother

\title{Saturating Constructions for Normed Spaces II}
\author{Stanislaw J. Szarek%
\thanks{Supported in part by a grant from the
National Science Foundation (U.S.A.).}
                      \and  Nicole Tomczak-Jaegermann%
\thanks{This author holds the Canada
         Research Chair in Geometric  Analysis.}
}

\newcommand\address{\noindent\leavevmode%
Department of Mathematics, 
Case Western Reserve University\\
Cleveland, Ohio 44106-7058, U.S.A.\\
{\it and}\\
Equipe d'Analyse Fonctionnelle, BP 186,
Universit\'e Pierre et Marie Curie\\
75252 Paris, France\\
{\small\tt%
szarek@ccr.jussieu.fr}\\[.5cm]
\noindent
Department of Mathematical and Statistical Sciences, 
University of Alberta\\
Edmonton, Alberta, Canada T6G 2G1,\\
{\small\tt%
                nicole@ellpspace.math.ualberta.ca} }
\date{}

\newtheorem{fact}{Fact}[section]
\newtheorem{thm}[fact]{Theorem}
\newtheorem{prop}[fact]{Proposition}
\newtheorem{lemma}[fact]{Lemma}

\newbox\nrmbox
\setbox\nrmbox=\hbox{$\Vert$}
\def\nrmrule{\vrule height\ht\nrmbox depth1.2\dp\nrmbox}
\setbox\nrmbox=%
                \hbox{\kern0.15em\nrmrule\kern0.15em\nrmrule\kern0.15em%
                          \nrmrule\kern0.15em}
\newcommand{\Snorm}[1]%
                {\copy\nrmbox#1\copy\nrmbox\kern-0.03em%
                      \lower.4ex\hbox{}}

\newcommand{\qed}{\bigskip\hfill\(\Box\)}

\newcommand{\R}{\mathbb{R}}
\newcommand{\Rn}{\R^n}

\newcommand{\E}{\mathbb{E}}
\newcommand{\pp}{\mathbb{P}}

\newcommand{\al}{\alpha}

\newcommand{\ep}{\varepsilon}

\newcommand{\ra}{{\rightarrow}}

\newcommand{\spn}{\mathop{\rm span\,}}

\newcommand{\rank}{\mathop{\rm rank\,}}

\newcommand{\Id}{\mathop{\rm Id}}
\newcommand{\conv}{\mathop{\rm conv\,}}
\newcommand{\convp}{\mathop{\rm conv}\nolimits{\!_p}}

\newcommand{\tr}{\mathop{\rm tr\,}}

\newcommand{\esp}{{\tilde{E}}}
\newcommand{\kk}{{\tilde{K}}}
\newcommand{\dd}{{\tilde{D}}}

\newcommand{\wrt}{with respect to }

\newcommand{\OO}{\Xi}

\begin{document}

\maketitle

{\abstract{\small
We prove several results of the following type:
given finite dimensional normed space $V$ possessing certain
geometric property there exists another space $X$ having the same property
and such that (1)
$\log \dim X = O(\log \dim V)$ and (2) every subspace of $X$, whose
dimension is not ``too small," contains a further well-complemented
subspace nearly isometric to $V$.  This sheds new light on the structure
of large subspaces or quotients of normed spaces (resp., large sections
or linear images of convex bodies) and provides definitive solutions to
several problems stated in the 1980s by V. Milman.}}

\section{Introduction}
\label{introduction}

This paper continues the study of the {\em saturation phenomenon } that
was discovered in \cite{ST1} and of the effect it has on our
understanding of the structure of high-dimensional normed spaces
and convex bodies. In particular, we obtain here a dichotomy-type result
which offers a near definitive treatment of some aspects of the
phenomenon.  We sketch first some background ideas and
hint on the broader motivation explaining the interest in the subject.

Much of geometric functional analysis revolves around the study of
the family of subspaces (or, dually, of quotients) of a given Banach
space.
In the finite dimensional case this has a clear geometric
interpretation: a normed space is determined by its unit ball,
a centrally symmetric convex body, subspaces correspond to sections
of that body,  and quotients to projections (or, more generally, linear
images). Such considerations are very natural from the geometric or
linear-algebraic point of view, but they also have a bearing
on much more applied matters. For example, a convex set may represent
all possible states of a physical system, and its sections or images
may be related to approximation or encoding schemes, or to results of an
experiment performed on the system.
It is thus vital to know to what degree the structure of
the entire space (resp., the entire set) can be recovered from the
knowledge of its subspaces or quotients (resp., sections/images).
At the same time, one wants to detect some  possible regularities in the
structure of subspaces  which  might have not existed in the whole space.

A seminal result in this
direction  is  the 1961 Dvoretzky theorem, with the 1971 strengthening
due to Milman, which  says that every symmetric convex body of
large dimension  $n$ admits
central sections which are approximately ellipsoidal and whose dimension
$k$ is of order $\log{n}$ (the order that is, in general, optimal).
Another major result  was the discovery of  Milman \cite{M1}
from the mid 1980's that {\em every }
$n$-dimensional normed space admits a {\em subspace of a quotient }
which is ``nearly" Euclidean and whose dimension is $\ge \theta n$,
where $\theta \in (0,1)$ is arbitrary (with the exact meaning of
``nearly" depending only on $\theta$).  Moreover,
a byproduct of the approach from \cite{M1} was the fact that every
$n$-dimensional normed space admits
a ``proportional dimensional'' quotient of bounded {\em volume ratio},
a volumetric  characteristic of a body
closely related to cotype properties (we refer to \cite{MS1}, \cite{T}
and \cite{P} for definitions of these and other basic
notions and results that are relevant here).
This showed that one can get a very essential regularity in a
global invariant of a space by passing to a
quotient {\em or} a subspace of dimension, say, approximately $n/2$.
It was thus natural to ask whether similar statements may
be true for other related characteristics. This line of thinking
was exemplified in a series of problems
posed by Milman in his 1986 ICM Berkeley lecture \cite{M2}.

The paper \cite{ST1} elucidated this circle of ideas and, in
particular, answered some of the problems from \cite{M2}.  A special 
but archetypal case of the main theorem from \cite{ST1} showed the
existence of an
$n$-dimensional space $Y$ whose {\em every } subspace (resp.,
{\em every } quotient) of dimension $\ge n/2$ contains a further
$1$-complemented subspace isometric to a preassigned (but a priori
arbitrary)
$k$-dimensional space $V$, as long as $k$ is at most of order $\sqrt{n}$.
In a sense, $Y$ was {\em saturated} with copies of the $V$.
This led to the discovery of the following phenomenon: passing to large
subspaces or quotients {\em can not}, in general, erase
$k$-dimensional features of a space if $k$ is below certain threshold
value depending on the dimension of the initial space and the exact
meaning of ``large."  In the particular case stated above, i.e., 
that of ``proportional" subspaces  or quotients, the
threshold dimension was (at least) of order $\sqrt{n}$,
and ``impossibility to erase" meant that every 
such subspace (resp.,  quotient map) 
preserved a copy of the given $V$.


However, the methods presented in \cite{ST1} were not sufficient for a
definitive treatment of the issue at hand. For example, we prove in
the present paper that, for any $q>2$, there are spaces of cotype $q$
(of arbitrarily high dimension $n$, with uniform control of constants)
whose all, say, $n/2$-dimensional subspaces are poorly $K$-convex (or,
equivalently, contain rather large subspaces well-isomorphic to
finite-dimensional $\ell_1$-spaces).  This is in stark contrast to
the extremal case of $q=2$: as it has been known since mid 1970's,
every space of cotype 2 admits proportional subspaces which are nearly
Euclidean (which is of course
incomparably stronger than being $K$-convex).
%
%
By comparison, in \cite{ST1} a similar result was established only for
$q >4$.  This answered one of the questions of Milman, but still left
open a possibility that an intermediate hypothesis weaker than cotype
2 (such as cotype $q$ with $2 < q \le 4$) could force existence of
nice subspaces. Our present theorem closes this gap completely, and
has the character of a dichotomy: for $q = 2$ every space of cotype 2
admits proportional nearly Euclidean subspaces, while for any $ q > 2$
there exist spaces of cotype $q$ without large $K$-convex subspaces at
all.  It was important to clarify this point since hypothetical
intermediate threshold values of $q$ (namely, $q=4$) appeared in related
-- and still not completely explained -- contexts in the asymptotic
geometric analysis literature, cf.
\cite{bourg} (see also Proposition 27.5 in \cite{T}) or \cite{pisier}.

Another variation of the saturation phenomenon that is being
considered here addresses what has being referred to recently as
``global properties." It has been realized in the last few years (cf.
\cite{MS2}) that many {\em local } phenomena (i.e., referring to
subspaces or quotients of a normed space) have {\em global }
analogues, expressed in terms of the entire space.  For example, a
``proportional" quotient of a normed space corresponds to the
Minkowski sum of several rotations of its unit ball. Dually, a
``proportional" subspace corresponds to the intersection of several
rotations.  (Such results were already implicit, e.g., in \cite{K}.)
Here we prove a sample theorem in this direction concerning the
Minkowski sum of two rotations of a unit ball, which answers a query
directed to us by V. Milman.

We use  the probabilistic method, and employ the ``blueprint" for
constructing random spaces that was developed by  Gluskin in
\cite{G} (the reader is also referred to \cite{MT1} for a survey of other
results and methods in this direction). In their most
general outline,  our arguments  parallel those of \cite{ST1}.
However, there are substantial differences,
and the present
considerations are much more subtle than those of \cite{ST1}.  Moreover,
we believe that several ingredients  (such as
a usage of Lemma~\ref{decoupling_lemma}-like statement to enable
decoupling of otherwise dependent events, or Lemma~\ref{2meanMstar}),
while playing mostly technical  role in this paper, are sufficiently
fundamental to be of independent interest.

The organization of the paper is as follows.  In the next section we
describe our main results and their immediate consequences.  We also
explain there the needed conventions employed by experts in the field, but
not necessarily familiar to the more general mathematical reader.
(Otherwise, we use the standard notation of convexity and geometric
functional analysis  as can be found, e.g., in \cite{MS1}, \cite{P} or
\cite{T}.)  Section \ref{2subspaces} contains the proof of Theorem
\ref{2hmmp}, relevant to the dichotomy mentioned above and to Problems
1-3 from \cite{M2}. Section \ref{2global-sec} deals with the global
variant of the saturation phenomenon.

\smallskip\noindent {\small {{\em Acknowledgement } Most of this
        research was performed while the second named author visited
        Universit\'e Marne-la-Vall\'ee and Universit\'e Paris 6 in the
        spring of 2002 and in the spring of 2003, and while both authors
        were attending the Thematic Programme in Asymptotic Geometric
        Analysis at the Pacific Institute of the Mathematical Sciences in
        Vancouver in the summer 2002.  Thanks are due to these
        institutions for their support and hospitality.} }

\section{Description of  results}
\label{2results}

The first result we describe is a {\em subspace} saturation theorem.
The approach of \cite{ST1} makes it easy to implement a saturation
property for subspaces.  Indeed, the dual space $X^*$ of the space
constructed in \cite{ST1}, Theorem 2.1 has the property that, under
some assumptions on $m$, $k$ and $n:= \dim X^*$, every $m$ dimensional
{\em subspace} of $X^*$ contains a (1-complemented) subspace isometric
to $V$ (where $V$ is a preassigned $k$-dimensional space).
In this paper we show that the construction can be performed while
preserving geometric features of the space $V$ (specifically, cotype
properties), a trait which is crucial to applications.


\begin{thm}
\label{2hmmp}
Let  $ q \in (2, \infty]$ and let $\ep >0$.  Then there exist
$\alpha = \alpha_q \in  (0, 1)$ and
$ c=c_{q, \ep}  > 0$
such that whenever positive integers $n$ and $m_0$ verify
$c^{-1}\,n^{\alpha} \le m_0 \le n$ and $V$ is {\bf any} normed space with
$$
\dim  V \le c m_0/ n^{\alpha} ,
$$
then there exists an $n$-dimensional nor\-med space $Y$ whose cotype
$q$ constant is bounded by a function of $q$ and the cotype $q$ constant
of $V$ and such that, for any \ $m_0 \le m \le n$,
{\bf every} $m$-dimensional subspace $\tilde{Y}$ of $Y$ contains a
$(1+\ep)$-complemented subspace $(1+\ep)$-isomorphic~to~$V$.
\end{thm}


Let us start
with several remarks concerning the hypotheses on $k := \dim V $ and
$m_0$ included in the statement above.  If, say, $m_0 \approx n/2$,
then $k$ of order ``almost'' $n^{1-\alpha}$ is allowed.  Nontrivial (i.e.,
large) values of $k$ are obtained whenever $m_0 \gg n^{\alpha}$;
we included the lower bound on $m_0$ in the statement to
indicate for which values of the parameters the assertion of the
Theorem is meaningful.

We can now comment on the relevance of Theorem \ref{2hmmp} to problems
from \cite{M2}.  Roughly speaking, Problems 2 and 3 asked whether
every space of nontrivial cotype $q < \infty$ contains a proportional
subspace of type $2$, or even just $K$-convex.  This is well known to be
true if
$q=2$ due to presence of nearly Euclidean subspaces
[For a reader not familiar with the type/cotype theory it will
be ``almost" sufficient to know that a nontrivial (i.e., finite) cotype
property of a space is equivalent to the absence of large subspaces
well-isomorphic to $\ell_\infty$-spaces; similarly,  nontrivial type
properties and $K$-convexity are related to the absence of
$\ell_1$-subspaces.]  Accordingly, by choosing, for example, $V =
\ell_1^k$ in the Theorem, we obtain -- in view of the remarks in the
preceding paragraph on the  allowed values of $k$ and $m$ -- a
space whose all ``large" subspaces  contain isometrically $\ell_1^k$
and which consequently provides a counter\-example to the problems for
{\em any } $q >2$.  More precisely,
if $m_0$ is ``proportional'' to $n$ and $V =
\ell_1^k$ is of the maximal dimension that is allowed, then the type 2
constant of any corresponding subspace $\tilde{Y}$ of $Y$ from the
Theorem is at least of order $n^{(1-\alpha)/2}$
(and analogously for any nontrivial type $p > 1$).
The $K$-convexity constant of any such $\tilde{Y}$ is
at least of order
$\sqrt{\log{n}}$ (up to a constant depending on $q$).
Problems 2 and 3 from
\cite{M2} are thus answered in the negative in a very strong sense.
Problem 1 from \cite{M2} corresponds to $q=\infty$ in Theorem \ref{2hmmp}
(i.e., no cotype assumptions) and has already been satisfactorily treated
in \cite{ST1}; however,  the present paper offers a unified
discussion of all the issues involved (see also related comments later in
this section).

We also remark that choosing $V = \ell_p^k$ (for some $1 < p <2$) in
Theorem \ref{2hmmp} leads to a space $Y$ whose type $p$ and cotype
$q$ constants are bounded by numerical constants and such that,
for every $m$-dimensional subspace $\tilde{Y}$ of $Y$ and every
$p < p_1 < 2$, the type $p_1$ constant of $\tilde{Y}$ is at least
$k^{1/p - 1/p_1}$.
If $m_0$ is
``proportional" to $n$, the type $p_1$ constant of $\tilde{Y}$ is at
least of order $n^{(1-\alpha)(1/p -   1/p_1)}$, in particular it tends to
$+\infty$  as $n \to \infty$. On the other hand, the spaces
$\tilde{Y}$ and $Y$ are then, by construction, uniformly (in
$n$) $K$-convex.

\smallskip


Theorem~\ref{2hmmp} will be an immediate consequence of the more
precise and more technical Proposition~\ref{2hmmptech} stated in
the next section. That statement makes the dependence
of the parameters $c$, $\alpha$ on  $\ep >0$ and $q \in (2,\infty)$
more explicit.
This will allow us, by letting $q \to \infty$, to retrieve the case
$q=\infty$ and then,
by passing to dual spaces, to reconstruct (up to a
logarithmic factor) the main theorem from \cite{ST1}: {\em if $n$,
$m_0$ and $k$ satisfy $\sqrt{n\, \log n} \le m_0 \le n$ and $k \le
m_0/\sqrt{n\, \log n}$, then for every $k$-dimensional normed space
$W$ there exists an $n$-dimensional normed space $X$ such that every
quotient $\tilde{X}$ of $X$ with $\dim \tilde{X} \ge m_0$ contains a
1-complemented subspace isometric to $W$.}


\smallskip

We wish now to offer a few comments on the construction that is
behind  Theorem \ref{2hmmp},
and which is implicit in  Proposition~\ref{2hmmptech}.
To this end, we recall some notation
and sketch certain ideas from \cite{ST1}, which also underlie the
present argument.

If $W$ is a normed space and $1 \le p < \infty$, by $\ell_p^N(W)$ we
denote the $\ell_p$-sum of $N$ copies of $W$, that is, the space of
$N$-tuples $(x_1, \ldots, x_N)$ with $x_i \in W$ for $1 \le i \le N$,
with the norm $\|(x_1, \ldots, x_N)\| = \left(\sum_i
    \|x_i\|^p\right)^{1/p}$. It is a fundamental and well-known
fact that the spaces $\ell_p^N(W)$ inherit type and cotype properties
of the space $W$, in the appropriate ranges of $p$ (cf. e.g., \cite{T},
\S 4).

The saturating construction from \cite{ST1} obtained $X^*$ as a
(random) subspace of $\ell_\infty^N(V)$, for appropriate value of $N$.
This is not the right course  of action in the context of Theorem
\ref{2hmmp} since such a subspace will typically contain rather
large subspaces  well-isomorphic to $\ell_\infty^s$, hence failing to
possess any nontrivial cotype property. However, substituting $q$ for
$\infty$ works:  the space
$\ell_q^N(V)$ and all its subspaces will be of cotype $q$ if $V$ is. The
approach of
\cite{ST1} was to concentrate on the case of $\ell_\infty^N(V)$, and then
to use the available ``margin of error" to transfer the results to  $q$
sufficiently close to $\infty$. By contrast, to handle the entire range
$2<q<\infty$ we need to work {\em directly } in the $\ell_q$ setting,
which -- as is well known to analysts -- often requires much more subtle
considerations.

\medskip

To state the next theorem, it will be helpful to subscribe to the
following ``philosophy" and notational conventions.
Since a normed space $X$ is completely described by its
unit ball $K=B_X$ or its norm $\|\cdot\|_X$, we shall tend to
identify these three objects.  In particular, we will write
$\|\cdot\|_K$ for the Minkowski functional defined by a
centrally symmetric  convex body $K \subset \R^n$ and denote the resulting
normed space  by $(\R^n, \|\cdot\|_K )$ or just $(\R^n, K )$.
Two normed spaces are isometric iff the corresponding convex bodies are
affinely equivalent.

\medskip

As suggested in the Introduction, it is of interest to consider
``global'' analogues of Theorem~\ref{2hmmp}-like statements.
  The following is a sample result that corresponds to the ``local"
Theorem 2.1 of \cite{ST1}, and that was already  announced in that paper.

\begin{thm}
\label{2hmm-global}
There exists a constant $c>0$ such that, for any positive integers $n,
k$ satisfying $k \le c n^{1/4}$ and
for every $k$-dimensional normed space $W$, there
exists an $n$-dimensional normed space $X= (\R^n, K)$ such that, for
any $u \in O(n)$, the normed space $(\R^n, K + u(K))$
contains a $3$-complemented subspace $3$-isomorphic  to $W$.
\end{thm}

In general, the interplay between the global and local results is not
fully understood.  While it is an experimental fact that a parallel
between the two settings exists, there is no formal conceptual
framework which explains it. It is thus important to provide more
examples in hope of clarifying the connection. It is also an
experimental fact that the local results and their global analogues
sometimes vary in difficulty. In the present context, the proof of
Theorem \ref{2hmm-global} is substantially more involved than that of
its local counterpart, Theorem 2.1 from \cite{ST1}.

\medskip


We conclude this section with several comments about notation.  As
mentioned earlier, our terminology is standard in the field and all
unexplained concepts and notation can be found, e.g., in \cite{MS1},
\cite{P} or \cite{T}.  The standard Euclidean norm on $\R^n$ will be
always denoted by $|\cdot|$.  ({\em Attention: } the same notation may
mean elsewhere cardinality of a set and, of course, the absolute value
of a scalar.) We will write $B_2^n$ for the unit ball in $\ell_2^n$
and, similarly but less frequently, $B_p^n$ for the unit ball in
$\ell_p^n$, $1 \le p \le \infty$.

For a set $S \subset \R^n$, by $\conv (S)$ we denote the convex hull
of $S$. If $1 \le p < \infty$, we denote by $\convp (S)$ the
$p$-convex hull of $S$, that is, the set of vectors of the form
$\sum_i t_i x_i$, where $t_i >0$ and $x_i \in S$ for all $i$, and
$\sum_i t_i^p =1$. (In particular, for $p=1$, $\convp (S) = \conv (S)$.)

The arguments below will use various subsets of $\R^n$ obtained as
convex hulls or $p$-convex hulls, for $1 < p < \infty$, of some more
elementary sets, or linear images of those; indeed for
Theorem~\ref{2hmmp} we have to consider the case of $p >1$, while in
Theorem~\ref{2hmm-global} the case of $p=1$ is sufficient.  In order
to emphasize the parallel roles which these sets (and other objects)
play in the proofs (which is also closely related to the role they
play in \cite{ST1}),  we try to keep a fully analogous notation for
them, and to distinguish them by adding a subscript $\cdot_p$ when the
set depends on $p$.

\section{Saturating  spaces of cotype $q >2$}
\label{2subspaces}

\noindent

Theorem \ref{2hmmp} will be an immediate consequence of the following
technical proposition.

\begin{prop}
\label{2hmmptech}
Let $2 < q < \infty$ and
set $\alpha := (q-2)/(2q+2)$ ($\in (0, 1/2)$).
Let $n$ and $m_0$ be positive integers with
$\sqrt q \, n^{1- \alpha} (\log{n})^{ (1-2\alpha)/3} \le m_0 \le n$.
Let $V$ be any normed space with
$$
\dim  V \le \frac{ c_1 m_0}{q^{1/2}\ n^{1- \alpha} (\log{n})^{
(1-2\alpha)/3}}
$$
(where $c_1 >0$  is an appropriate  universal constant).
Then there exists an
$n$-dimensional nor\-med space $Y$ whose cotype
$q$ constant is bounded by a function of $q$ and the cotype $q$ constant
of $V$ and such that, for any \ $m_0 \le m \le n$,
every $m$-dimensional subspace $\tilde{Y}$ of $Y$ contains a
$2^{1/q}$-complemented subspace $2^{1/q}$-isomorphic~to~$V$.
Moreover, for every $\ep>0$, we may
replace the quantity $2^{1/q}$ 
by $1+\ep$,
at the cost of allowing $c_1$  to depend on $\ep$.
\end{prop}

\noindent
{\em Proof }
Fix $2 < q < \infty$ and let $p = q/(q-1)$ be the
conjugate exponent. Let $1\le k \le m \le n \le kN$ be
positive integers.  More restrictions will be added on these
parameters as we proceed, and in particular we shall specify $N$
(depending also on $q$) at the end of the proof.  Notice that choosing
the constant $c_1$ small makes the assertion vacuously satisfied for
small values of $m_0$, and so we may and shall assume that $m_0, n$
and $N$ are large.

Let $V$ be a $k$-dimensional normed space. Identify $V$ with $\R^k$ in
such a way that the Euclidean ball $B_2^k$ and the unit ball $B_V$ of
$V$ satisfy $ B_2^k \subset B_V \subset \sqrt k \, B_2^k$ (for example,
$B_2^k$ may be the ellipsoid of maximal volume contained in $B_V$).
As indicated in the preceding section, we shall construct the space
$Y$ as a (random) subspace of $\ell_q^N (V)$, the $\ell_q$-sum of $N$
copies of $V$.  We will actually work in the dual setting of random
quotients of $Z_p := \ell_p^N (W)$, where $W:=V^*$; as frequent in this
type of constructions, the geometry of that setting is more transparent.
The above identification of $V$ with $\R^k$ induces the identification of
$W$ with $\R^k$, and thus allows to identify $Z_p$ with $\R^{Nk}$.

Let $G = G(\omega)$ be a $n \times Nk$ random matrix (defined on some
underlying probability space $(\Omega, \pp)$) with independent $N(0,
1/n)$-distributed Gaussian entries. Consider $G$ as a linear operator
$G: \R^{Nk} \to \R^n$ and set
\begin{equation}
\label{2def_body}
K_p= B_{X_p(\omega)}  := G (\omega) (B_{Z_p}) \subset \R^n.
\end{equation}
The random normed space $X_p = X_p(\omega)$ can be  thought of
as a random (Gaussian) quotient of $Z_p$, with $G (\omega)$ the
corresponding quotient map and $K_p$ the unit ball of $X_p$.  [The
normalization of
$G$ is not important; here we choose it so that, with $k, N$ in the ranges
that matter, the  radius of the Euclidean ball circumscribed on  $K_p$
be typically comparable to 1.]

We reiterate that the dual spaces $X_p^* = X_p(\omega)^*$ are
isometric to subspaces of $Z_p^* = \ell_q^N (V)$ and so their cotype
$q$ constants are uniformly bounded (depending on $q$ and the cotype
$q$ constant of $V$).
We shall show that, for appropriate choices of the parameters, the
space $Y= X_p(\omega)^*$ satisfies, with probability close to $1$, the
(remaining) assertion of Theorem \ref{2hmmp} involving the subspaces
well-isomorphic to $V$.  This will follow if we show that, outside of a
small exceptional set, every quotient $\tilde{X}_p(\omega)$ of
$X_p(\omega)$ of dimension $m \ge m_0$ contains a
$2^{1/q}$-complemented subspace $2^{1/q}$-isomorphic to $W$, for
values of $k$ described in Theorem \ref{2hmmp} (and analogously for
$1+\ep$ in place of $2^{1/q}$). To be absolutely precise, we shall show
that the identity on $W$ well factors  through $X_p(\omega)$, a
property which dualizes without any loss of the constant involved.
Thus we have a very similar problem to the one considered in
\cite{ST1}, however the present context requires several subtle
technical modifications of the argument applied there.

Similarly as in \cite{ST1}, we will follow the scheme first employed
in \cite{G}: {\em Step~I\ } showing that the assertion of the theorem
is satisfied for a {\em fixed } quotient map with probability close to 1;
{\em Step II\ } showing that the assertion is ``essentially stable" under
small perturbations of the quotient map; and {\em Step III\ } which
involves a discretization argument.


\medskip

We start by introducing some notation that will be used throughout the
paper.  Denote by $F_1, \ldots, F_N$ the $k$-dimensional coordinate
subspaces of $\R^{Nk}$ corresponding to the consecutive copies of
$W$ in $Z_p$.  In particular, from the definition of the $\ell_p$-sum
we have
$$
B_{Z_p} = \convp (F_j \cap B_{Z_p} : j \in \{1, \ldots , N \}).
$$
For $j=1,
\ldots, N$, we define subsets of $\R^n$ as follows: $E_j := {G}(F_j)$,
$K_j := {G}(F_j \cap B_{Z_p})$ and
\begin{equation} \label{2K'}
K_{j,p}' := {G}(\spn[F_i: i \ne j] \cap B_{Z_p}).
= \convp (K_i: i \ne j).
\end{equation}
We point out certain ambiguity in the notation: $K_p$, $p\in (1,2)$,
is the unit ball of $X_p$, while $K_j$, $j\in \{1, \ldots , N \}$
stands for the section of $K_p$ corresponding to $E_j$. This
should not lead to confusion since, first, the sections {\em do not}
depend on $p$ and, second, $p$ remains fixed throughout the argument.
(Similar caveats apply to the families of sets $D_\cdot$,
$\kk_\cdot$ and $\dd_\cdot$ which are defined in what follows.)

In addition to $K_p$ and the $K_j$'s, we shall need subsets
constructed in an analogous way from the Euclidean balls. First, for
$j=1, \ldots, N$, set
$D_j := {G}(F_j \cap  B_2^{Nk})$. Then let
\begin{equation}
   \label{dp}
D_p := G \left( \convp
(F_j \cap B_2^{Nk} : j \in \{1, \ldots , N \})\right)
=  \convp (D_j: j \in \{1, \ldots , N \}).
\end{equation}
Next, for $j=1, \ldots, N$, let
\begin{equation}
   \label{djp'}
     D_{j, p}'
:= G \left( \convp
(F_i \cap B_2^{Nk} :  i \ne j )\right)
=  \convp (D_i: i \ne j).
\end{equation}
Finally,  for a subset $I \subset  \{1, \ldots, N\}$, we let
\begin{equation}
\label{dIp}
D_{I, p}
:= G \left( \convp
(F_i \cap B_2^{Nk} :  i \in I )\right)
=  \convp (D_i: i \in I).
\end{equation}

Note that since $\frac1{\sqrt k} B_2^k \subset B_{W} \subset B_2^k$,
it follows that $\frac 1 {\sqrt{k}}D_j \subset K_j \subset D_j$.
Consequently, analogous inclusions hold for all the corresponding $K$-
and $D$-type sets as they are $p$-convex hulls of the appropriate
$K_j$'s and $D_j$'s.

\bigskip
\noindent {\em Step I.\ Analysis of a single quotient map. }  Since
a quotient space is determined up to an isometry by the kernel of a
quotient map, it is enough to consider quotient maps which are {\em
orthogonal } projections.  Let, for the time being, $Q: \R^n \to
\R^m$ be the canonical projection on the first $m$ coordinates.  In view
of symmetries of our probabilistic model, all relevant features of this
special case will transfer to an arbitrary rank $m$ orthogonal projection.

Let $\tilde{G} = Q G$, i.e.,  $\tilde{G}$ is
the $m \times Nk$ Gaussian matrix obtained by restricting $G$ to
the first $m$ rows.
Let $\kk_p = Q(K_p) = \tilde{G}(B_{Z_p})$ and denote
the space $(\R^m, \kk_p)$ by $\tilde{X}_p$; the space $\tilde{X}_p$
is the quotient of $X_p$ induced by the quotient map $Q$.
We shall use the notation of
$\esp_j$,  $\kk_j$,   $\kk_{j,p}'$ 
for the subsets of $\R^m$ defined in the same way
as $E_j$, $K_j$, $K_{j,p}'$, 
above, but using the matrix $\tilde{G}$ in place of $G$.
Analogous convention is used to define the
$\dd$-type sets
$\dd_p$,  $\dd_j$ and    $\dd_{j,p}'$.

\smallskip

For any subspace $H\subset \R^m$, we will denote by $P_H$
the orthogonal projection onto $H$.
We shall show that outside of an exceptional set of small measure there
exists
$j \in \{1, \ldots, N\}$ such that $P_{\esp_j} ( \kk_{j,p}') \subset
\kk_j$.
Note that, for any given $i$, we {\em always } have $\kk_p = \convp
(\kk_i, \kk_{i,p}')$
and $\kk_i \subset \esp_i$. It follows that, for $j$
as above,
\begin{equation}
P_{\esp_j}( \kk_p) = \convp (\kk_j, P_{\esp_j}( \kk_{j,p}')) \subset
2^{1/q}    \kk_j.
\label{2koniec}
\end{equation}
Note that $\kk_j$ is an affine image of the ball $F_j \cap
B_{Z_p}$, which is the ball $B_W$ on coordinates from $F_j$.  On
the other hand, $\esp_j $ considered as a subspace of $\tilde{X}_p$
(thus endowed with the ball $\esp_j \cap \kk_{p}$) satisfies, by
(\ref{2koniec}), $ \kk_j \subset \esp_j \cap \kk_{p} \subset
2^{1/q} \kk_j$, which makes it $2^{1/q}$-isomorphic to $B_W$.
Using (\ref{2koniec}) again we also get the $2^{1/q}$-complementation.
(Similarly, $P_{\esp_j} ( \kk_{j,p}') \subset \ep \kk_j$ will imply
$(1+\ep)$-isomorphism and $(1+\ep)$-complementation.)

Returning to inclusions between the $K$- and $D$-type sets, they also
hold for the
$\kk$- and $\dd$-type sets,
so that, for example, $\frac 1
{\sqrt{k}}\dd_j \subset \kk_j \subset \dd_j$.
Consequently,
in order for the inclusion
$P_{\esp_j} (\kk_{j,p}') \subset \kk_j$ to hold it is enough to have
\begin{equation}
\label{2brutal}
P_{\esp_j}( \dd_{j,p}') \subset \frac 1 {\sqrt{k}}\dd_j.
\end{equation}
The rest of the proof of Step I is to show that, with an appropriate
choice of the parameters, this seemingly rough condition is satisfied
for some $1 \le j \le N$, outside of a small exceptional set.

\medskip

Let us now pass
to the definition of the exceptional set.  We start by
introducing, for $j \in \{1, \ldots, N\}$, the ``good'' sets.  Fix a
parameter $ 0 < \kappa \le 1$ to be determined later, and let
\begin{eqnarray}
\label{Theta'j}
\Theta_{j}' & := &
\left\{ \omega \in \Omega :
P_{\esp_j} (\dd_{j,p}') \subset
\kappa B_2^m \right\} \\
\Theta_{j,0}' &:= & \left\{ \frac12 \sqrt{\frac m{n}}(B_2^m \cap
\esp_j) \subset \dd_j
\subset
2\sqrt{\frac m{n}}(B_2^m \cap \esp_j) \right\}.
\label{Theta''(j)}
\end{eqnarray}

Now if  $\kappa $, $k$,  $m$ and $n$ satisfy
\begin{equation} \label{2cubeball}
\kappa \leq \frac 1
      {\sqrt{k}} \cdot \frac12 \sqrt{\frac
m{n}} \, ,
\end{equation}
then, for
$\omega \in \Theta_{j}' \cap \Theta_{j,0}' $,
the inclusion
(\ref{2brutal}) holds.
Thus, outside of the  exceptional set
\begin{equation}
\label{Theta0}
\Theta^0 :  =  \Omega \setminus \bigcup_{1\le j \le N}
\bigl(\Theta_{j}' \cap \Theta_{j, 0}'\bigr)
= \bigcap_{1\le j \le N} \left( (\Omega \setminus \Theta_{j}')
\cup
(\Omega \setminus \Theta_{j,0}')\right)
\end{equation}
there exists $j \in \{1, \ldots, N \}$ such that (\ref{2brutal})
holds, and this implies, by an earlier argument, that there exists $j
\in \{1, \ldots, N \}$ such that $\esp_j $ considered as a subspace of
$\tilde{X}_p$ is $2^{1/q}$-complemented and $2^{1/q}$-isomorphic to
$W$.

\smallskip

It remains to show that the measure of the exceptional set $\Theta^0 $
is appropriately small; this will be the most technical part of the
argument.  The first problem we face is that the events entering the
definition of $\Theta^0 $ are not independent as $j$ varies. We overcome
this difficulty  by a decoupling trick
which allows to achieve conditional independence on a large subset of
these events.

\begin{lemma}
\label{decoupling_lemma}
Let $\Lambda = (\lambda_{ij})$ be an $N \times N$ matrix such that
\newline $1^\circ$ $0 \le \lambda_{ij} \le 1$ for all $i, j$
\newline $2^\circ$ $\sum_{i=1}^N \lambda_{ij} =1$ for all $j$
\newline $3^\circ$ $\lambda_{jj} =0$ for all $ j$.
\newline Then there exists $J \subset \{1, \ldots, N\}$ such that
$|J| \ge N/3$  and for every $j \in J$ we have
$$
\sum_{i \not\in J} \lambda_{ij} \ge 1/3.
$$
\end{lemma}

\noindent
{\em Proof } This lemma is an immediate consequence of the result of
K.~Ball on suppression of matrices presented and proved in \cite{BT}.
By Theorem 1.3 in \cite{BT} applied to $\Lambda$, there exists a
subset $J \subset \{1, \ldots, N\}$ with $|J| \ge N/3$ such that
$\sum_{i \in J} \lambda_{i j} < 2/3$ for $j \in J$, which is just a
restatement of the condition in the assertion of the Lemma.
\qed

Now if $\omega \in \Omega \setminus \Theta_{j}'$, for some $j = 1,
\ldots, N$, then, by  (\ref{Theta'j}) and
the definition of $\dd_{j,p}'$, there
exist  $x_{i,j} \in F_i\cap B_2^{Nk}$, for all $i \neq j$, with
$\sum_{i\neq j} |x_{i,j}|^p =1$ and $z_j \in \esp_j \cap B_2^m $ such that
$$
\langle \tilde{G} \bigl( \sum_{i \ne j} x_{i,j}\bigr), z_j
\rangle =:\kappa_j > \kappa.
$$
By changing $x_{i,j}$ to $- x_{i,j}$ if necessary, we may assume
that $(\tilde{G}x_{i,j}, z_j ) \ge 0$ for all $ i \ne j$.  Thus if
$\omega \in \Omega \setminus \Theta_{j}'$, for all $j$, then we can
consider the matrix $\Lambda$ defined, for $j=1, \ldots, N$, by
$\lambda_{ij} = \langle \tilde{G}x_{i,j}, z_j \rangle  /\kappa_j$ for $i
\ne j$ and
$\lambda_{jj} =0$.  Let $J \subset \{1, \ldots, N\}$ be the set
obtained by Lemma~\ref{decoupling_lemma}.  Then $|J| \ge N/3$ and for
every $j \in J$ we have
$$
\langle  \tilde{G} \bigl( \sum_{i \not\in J} x_{i.j}\bigr), z_j \rangle
\ge \kappa/3,
$$
and so
$$
\omega \in \Omega \setminus
\left\{ \omega \in \Omega :
P_{\esp_j} (\dd_{J^c, p}
)    \subset (\kappa/3) B_2^m \right\};
$$
we recall that  for a subset $I \subset  \{1, \ldots, N\}$,
$\dd_{I, p}$ has been defined in (\ref{dIp}),
and $J^c = \{1, \ldots, N\} \setminus J$.

Let ${\mathcal J}$ be the family of all subsets $J \subset
\{1, \ldots, N\}$ with $|J |= \lceil N/3 \rceil =:\ell$.
Then the above argument immediately
implies that
\begin{equation}
  \label{big-incl}
\bigcap_{1 \le j \le N} (\Omega \setminus \Theta_{j}')
\subset \bigcup_{J \in {\mathcal J}}
\bigcap_{j \in J} (\Omega \setminus \Theta_{j, J^c}'),
\end{equation}
where
\begin{equation}
  \label{ThetaJc}
\Theta_{j, J^c}' :=
\left\{ \omega \in \Omega :
P_{\esp_j} (\dd_{J^c, p}) \subset
(\kappa/3) B_2^m \right\}.
\end{equation}

This definition has a form similar to
(\ref{Theta'j}) (indeed, $\dd_{j,p}' = \dd_{I, p}$ where $I = \{1,
\ldots, N\} \setminus \{j\}$; additionally,  $\kappa$ gets replaced by
$\kappa/3$).
Comparing (\ref{big-incl}) with (\ref{Theta0}) and
reintroducing the sets $\Omega \setminus \Theta_{j, 0}'$  into our
formulae we obtain
\begin{equation}
  \label{Theta0cover}
  \Theta^0 \subset  \bigcup_{J \in {\mathcal J}}
\Theta_{J} ,
\end{equation}
where for $J \in {\mathcal J}$ we set
\begin{equation}
  \label{ThetaJ}
\Theta_J := \bigcap_{j \in J} (\Omega \setminus (\Theta_{j, J^c}'
\cap \Theta_{j, 0}')).
\end{equation}

Our next objective will be to estimate $\pp(\Theta_J)$ for a fixed $J$.
By symmetry, we may restrict our attention to
$J = \{1, \ldots, \ell \}$.
For  $j = 1, \ldots, \ell$, set
$$
{\mathcal E_{j, p}} :=
\Omega \setminus (\Theta_{j, J^c}'
\cap \Theta_{j, 0}').
$$
Then
$$
\Theta_J  = \bigcap_{j \in J} {\mathcal E_{j, p}}.
$$
We are now in the position to make the key observation of this part of
the argument: for a fixed $J \in {\mathcal J}$, the events
$\mathcal{E}_{j, p} $, for $j \in J$, are conditionally independent
with respect to $\dd_{J^c,p}$: once
$\dd_{J^c,p}$ is fixed, each $\mathcal{E}_{j,p} $ depends only
on the restriction $G_{|F_j}$.  In fact, the ensemble
$\{\tilde{G}_{|F_j} : j \in J \}
\cup \{\dd_{J^c,p}\}$ is independent since its distinct elements
depend on disjoint sets of columns of $\tilde{G}$,
and the columns themselves
are independent.  This and the symmetry in the indices $j \in J$
implies that
\begin{eqnarray}
\label{decoupling-p}
\pp(\Theta_J \mid \dd_{J^c,p} )
& = &\pp\Bigl(\bigcap_{j \in J} (\mathcal{E}_{j,p}
\mid \dd_{J^c,p} )\Bigr) \nonumber\\
&=& \prod_{j \in J} \pp(\mathcal{E}_{j,p}
  \mid \dd_{J^c,p} )
=\Bigl(\pp( \mathcal{E}_{1,p}  \mid \dd_{J^c,p} )\Bigr)^{\ell}.
\end{eqnarray}

\medskip

To estimate $\pp( \mathcal{E}_{1,p}  \mid \dd_{J^c,p} )$ first
note that, by the definition
of  $\mathcal{E}_{1,p}$ this probability is less than or equal to
$\pp(\Omega \setminus \Theta'_{1,J^c} \mid \dd_{J^c,p})
+ \pp(\Omega \setminus \Theta_{1,0}'  \mid \dd_{J^c,p} )$.
Next, since $\Theta_{1,0}'$ is independent of $\dd_{J^c,p}$, the second
term equals just $1 - \pp( \Theta_{1,0}')$. Further, the set
$\Theta_{1,0}'$ is the same as in \cite{ST1} (where it was denoted by
$\Omega_{1,0}'$, see formula (3.7) in that paper), and so
\begin{equation}
\label{2bigcube}
\pp(\Omega \setminus \Theta_{1,0}'  \mid \dd_{J^c,p} )
\le  e^{-m/32} + e^{-9m/32};
\end{equation}
(see (3.16) in \cite{ST1}, or use directly Lemma~3.3 from \cite{ST1}
or Theorem 2.13 from \cite{DS}, both of which describe the behaviour of
singular numbers of rectangular Gaussian matrices).

\smallskip For the term involving $\Omega \setminus
\Theta_{1, J^c}'$ the probability estimates are much more delicate and
will require two auxiliary lemmas. Before we state them, we
recall the by now classical concept of functional $M^*(\cdot)$,
defined for a set $S \subset \R^d$ by
\begin{equation}
   \label{2Mstar}
M^*(S) :=
\int_{S^{d-1}}{\sup_{y \in S}
{\langle x, y \rangle}} dx ,
\end{equation}
where the integration is performed with respect to the normalized
Lebesgue measure on $S^{d-1}$ (this is $1/2$ of what geometers call
the mean width of $S$; if $S$ is the unit ball for some norm, $M^*(S)$
is the average of the {\em dual } norm over $S^{d-1}$).
We then have

\begin{lemma}
\label{2meanMstar}
Let $d, s$ be integers with $1\le d \le s$ and let
$A = (a_{ij})$ be a
$d \times s$ random matrix with independent
$N(0, \sigma^2)$-distributed Gaussian entries.
Further, let
$a >0$ and let $S \subset \R^s$
be a symmetric convex body
satisfying  $S \subset a B_2^s$.
Then the random body \ $A S\subset \R^d$  verifies
$$
{\E}\left(M^*(A S)\right) = c_s \sigma M^*(S),
$$
where $c_s =\sqrt{2}\Gamma(\frac{s+1}2)/\Gamma(\frac s2)\le \sqrt{s}$.
Moreover, for any $t > 0$,
$$
\pp\left(M^*(A S) > c_s \sigma M^*(S) + t \right)
\le e^{-dt^2/2a^2\sigma^2}.
$$
\end{lemma}
{\em Proof \ }  The first assertion is quite standard. We have,
$$
{\E}\left(M^*(A S)\right) =
{\E}\int_{S^{d-1}}{\sup_{x \in S}
{\langle A x, y \rangle}} dy =
\int_{S^{d-1}}{\E}\,{\sup_{x \in S}
{\langle x, A^* y \rangle}} dy .
$$
Since, for any $y \in \R^d$,  $A^* y$  is distributed as $\sigma |y|$
times  the standard Gaussian vector in $\R^s$, the integrand
${\E}\,{\sup_{x \in S} {\langle x, A^* y \rangle}}$
does not depend on $y \in {S^{d-1}}$ and is equal to
the appropriate (independent of $S$) multiple of the spherical mean.
The value of the $c_s$ may be obtained, e.g., by calculating
the Gaussian average for $S=S^{s-1}$.

For the second assertion, we show first that the function $T \ra f(T)
:= M^*(T S)$ is $a/\sqrt{d}$-Lipschitz \wrt the Hilbert-Schmidt norm
$\|\cdot\|_{HS}$.  Indeed, directly from the definition (\ref{2Mstar})
we have
\begin{eqnarray*}
f(T_1) - f(T_2) &= &\int_{S^{d-1}}{\sup_{x \in S}
{\langle T_1 x, y \rangle}} dy - \int_{S^{d-1}}{\sup_{x \in S}
{\langle T_2 x, y \rangle}} dy \\
&\le  & \int_{S^{d-1}}{\sup_{x \in S}
{\langle (T_1-T_2) x, y \rangle}} dy \\
&\le  & \int_{S^{d-1}} a |(T_1-T_2)^*y| dy \\
&\le  & a \left( \int_{S^{d-1}} |(T_1-T_2)^*y|^2 dy \right)^{1/2}\\
&=& (a/\sqrt d) \|(T_1-T_2)^*\|_{HS} = (a/\sqrt d) \|T_1-T_2\|_{HS}.
\end{eqnarray*}
The Gaussian isoperimetric inequality
(see  e.g.,  \cite{L},  (2.35))
implies now that
$$
\pp\bigl(M^*(A S) > {\E}\left(M^*(A S)\right) + t \bigr)
\le e^{-dt^2/2a^2\sigma^2},
$$
which shows that the second assertion of the Lemma follows
from the first one.  \qed

The second lemma describes the behaviour of the diameter
of a random rank $d$ projection (or the image under a Gaussian map)
of a subset of
$\R^s$. Let $d, s$ be integers with  $1 \le d \leq s$ and
let $G_{s,d}$ be the Grassmann manifold of $d$-dimensional
subspaces of $\R^s$ endowed with the normalized Haar measure.

\begin{lemma}
\label{shrinking}
Let  $a >0$ and let $S \subset \R^s$ verify \ $S \subset a B_2^s$.
Then, for any  $t > 0$, the set
$\left\{ H \in G_{s,d} : P_H (S) \subset
\left(a \sqrt{d/s} + M^*(S) + t \right) B_2^s \right\}$
has measure $\geq 1- \exp(-t^2s/2a^2 + 1)$. Similarly, replacing
$P_H$ by a $d \times s$ Gaussian matrix $A$ with independent $N(0,1/s)$
entries, we obtain a lower bound on probability of the form
$1- \exp(-t^2s/2a^2)$.
\end{lemma}
The phenomenon
discussed in the Lemma is quite well known, at least if
we do not care about the specific values of numerical constants (which
are not essential for our argument) and precise estimates on
probabilities.  It is sometimes refereed to as the ``standard
shrinking'' of the diameter of a set, and it is implicit, for example,
in probabilistic proofs of the Dvoretzky theorem, see \cite{M-dv},
\cite{MS1}. A more explicit statement can be found in \cite{M3},
and the present version was proved in \cite{ST1}.

\medskip
We now return to the main line of our argument.
Observe first that for $1 < p \le 2$ one has
\begin{equation}
\label{cp_const}
M^* \left( \convp (F_j \cap B_2^{Nk}: 1 \le j \le N)\right)
= M^* \bigl( \ell_p^N (\ell_2^{k})\bigr)
\le  C_p N^{1/q - 1/2}, 
\end{equation}
where $q = p/(p-1)$ and $1 \le C_p \le C \sqrt q$, where $C >0$ is an
absolute constant. This is likely known, and certainly follows by
standard calculations; e.g., by passing to the average of the
$\ell_q^N (\ell_2^{k})$-norm (dual to the $\ell_p^N
(\ell_2^{k})$-norm; cf.  the comments following (\ref{2Mstar})),
expressing it in terms of the Gaussian average and then majorizing the
latter via the $q$th moment, which in turn may be explicitly computed.

\smallskip The estimate (\ref{cp_const}) has two consequences for the
set $D_p$ (defined in (\ref{dp})). Firstly, the Gaussian part of Lemma
\ref{shrinking} implies that, with our normalization of $G$, the diameter
of $D_p$ is typically comparable to 1. More precisely, consider the
exceptional set
\begin{equation}
\label{Theta1}
  {\Theta}^1 := \{ \omega : D_p \not\subset 2 B_2^n\}.
\end{equation}
Then, as long as $M^* \bigl( \ell_p^N
(\ell_2^{k})\bigr) \le (n/(4Nk))^{1/2}$,
we can apply Lemma \ref{shrinking} to the
$n \times Nk$ \ matrix $A = (n/(Nk))^{1/2}G$ and $t = 1/2$
to obtain $\pp({\Theta}^1) \le \exp (-n/8)$
(note that $a =1$ in this case).  On the other hand, by (\ref{cp_const}),
the estimate on  $M^* \bigl( \ell_p^N(\ell_2^{k})\bigr)$ is
satisfied whenever
\begin{equation} \label{fortheta1}
n \ge 4 C_p^2 N^{2/q}k ,
\end{equation}
which will be ensured by our final choice of $N$ and the conditions
that will be imposed on the dimensions involved.

Secondly, by Lemma~\ref{2meanMstar}, we have ${\E}\left (M^*
(D_p)\right) \le C_p \sqrt{k/n}\, N^{1/q}$.  (Recall that $G$ is an $n
\times Nk$ Gaussian matrix with $\sigma^2 = 1/n$.)  Thus, by
the second part of the Lemma, our second exceptional set
\begin{equation}
\label{2Thetabar}
  \bar{\Theta}^1 := \{ \omega : M^*(D_p) >
              2  C_p \sqrt{k/n}\, N^{1/q}\}
\end{equation}
satisfies $\pp (\bar{\Theta}^1) \le  \exp(- C_p^2 k n N^{2/q}/2) \le
\exp (-n/2)$  (remember that $C_p \ge 1$).

\medskip

Now recall that $Q: \R^n \to \R^m$ is the canonical projection on the
first $m$ coordinates.
Since $\dd_p = Q D_p$, it follows that for $\omega \not\in \Theta^1$ we
have
\begin{equation}
\label{2ddtildep}
\dd_{J^c,p} \subset \dd_p \subset 2 B_2^m.
\end{equation}
Further, it is a general fact (shown by passing to Gaussian averages)
that, for any $S \subset \R^n$,
$ c_m M^* (Q S) \le c_n M^* (S)$, were $c_m$ and $c_n$ are constants
from Lemma~\ref{2meanMstar} and $c_n / c_m \le (2/\sqrt \pi)\sqrt
{n/m}$. Thus, for $\omega \not\in \bar{\Theta}^1$ we have
\begin{equation}
\label{2ddtildepMstar}
M^*(\dd_{J^c,p} ) \le
M^* (\dd_p) \le  \frac{2}{\sqrt \pi}\sqrt {\frac{n}{m}}\cdot 2 C_p
\sqrt {\frac{k}{n}}
\,N^{1/q}
     = C_p \frac{4}{\sqrt \pi} \sqrt {\frac{k}{m}}\, N^{1/q}.
\end{equation}

We now return to our current main task, which is to analyze
the set $\Omega \setminus \Theta_{1, J^c}'$.  Since we are working with
conditional probabilities, we need to introduce another exceptional set
which is $\dd_{J^c, p}$-measurable
\begin{equation}\label{Theta'}
\Theta' := \{ \omega : \dd_{J^c, p} \not \subset 2 B_2^m \mbox{ or }
M^* (\dd_{J^c, p}) > C_p (4/\sqrt \pi) \sqrt {k/m}\, N^{1/q} \}
\end{equation}
It follows directly from
(\ref{2ddtildep}) and (\ref{2ddtildepMstar}) that $\Theta' \subset
{\Theta}^1 \cup\bar{\Theta}^1$. We emphasize that $\Theta'$ depends in
fact on $J$, but $J$ is fixed at this stage of the argument. Moreover, the
sets $\Theta'$ corresponding to different $J$'s are subsets of a small
common superset ${\Theta}^1 \cup\bar{\Theta}^1$, which is additionally
independent of $Q$.

\smallskip The definition of the set $\Theta_{1, J^c}'$
(cf. (\ref{ThetaJc})) involves the
diameter of a random rank $k$ projection of $\dd_{J^c, p}$
(note that, by the rotational invariance of the Gaussian measure,
${\esp_1}$ is distributed uniformly in $G_{m,k}$ , and is independent
of $\dd_{J^c, p}$).
Moreover, if  $\omega \not\in \Theta'$,
we control the diameter and $M^*$ of the set $S = \dd_{J^c, p}$, and so we
are  exactly in a position to apply Lemma \ref{shrinking}.
Specifically, we use $t = \kappa/6$, $a=2$  and assume that
\begin{equation}
\label{lowerkappa}
C_p (4/\sqrt \pi)
\sqrt{k/m}\, N^{1/q} \le \kappa/12
\end{equation}
(which implies $a \sqrt{k/m} = 2 \sqrt{k/m} \le \kappa/12$) to obtain
\begin{equation}
\label{2Theta_est}
\pp  \left( \Omega \setminus \Theta_{1, J^c}' | \dd_{J^c, p} \right) \le
e^{- \kappa^2 m / (8 \cdot 6^2) +1}.
\end{equation}

For the record, we note that  (\ref{lowerkappa}) implies
$$
4 C_p^2 N^{2/q}k \le 4 (\kappa/12)^2 m \le m \le n,
$$
and thus the condition (\ref{fortheta1}) that appeared in connection with
the measure estimate for $\Theta^1$ is automatically satisfied.

\smallskip

Substituting (\ref{2Theta_est}) combined with the estimate (\ref{2bigcube})
for the measure of $\Omega \setminus \Theta_{1, 0}'$ into
(\ref{decoupling-p}) we deduce that, outside of ${\Theta}'$,
\begin{eqnarray*}
\pp(\Theta_J \mid D_{J^c, p} )
& \le  & \left(e^{-\kappa^2 m/(8 \cdot 6^2) +1}+ e^{-m/32} +
e^{-9m/32}\right)^{\ell}\\
& \le & (2e)^{\ell} e^{-\kappa^2 m \ell/(8 \cdot 6^2)}.
\end{eqnarray*}
Averaging over $\Omega \setminus {\Theta}'$
(and using $\Theta' \subset {\Theta}^1 \cup\bar{\Theta}^1$) yields
$$
\pp\left(\Theta_J \setminus ({\Theta}^1 \cup
\bar{\Theta}^1)\right)
\le \pp\left(\Theta_J \setminus {\Theta}'\right)
\le \pp\left(\Theta_J \mid \Omega \setminus {\Theta}'\right)
\le (2e)^{\ell} e^{-\kappa^2 m \ell/(8 \cdot 6^2)}.
$$
%
Since
$|\mathcal{J}| = {{N}\choose{\ell}}$ and $\bigcup_{J \in \mathcal{J}}
\Theta_J \supset \Theta^0$ (cf. (\ref{Theta0cover})),  it follows that
\begin{equation} \label{1singleop-p}
\pp\left( \Theta^0 \setminus
({\Theta}^1 \cup \bar{\Theta}^1)
\right) \le
\pp\left(\bigcup_{J \in \mathcal{J}} \Theta_J \setminus
({\Theta}^1 \cup \bar{\Theta}^1)
\right) \le {{N}\choose{\ell}}  (2e)^\ell
e^{-\kappa^2 m \ell/(8 \cdot 6^2)}.
\end{equation}
Consequently ,
\begin{eqnarray}
\label{2singleop-p}
\pp(\Theta^0) &\le&  \pp(\Theta^1) +
\pp(\bar{\Theta}^1) +
\pp\left(\Theta^0 \setminus
({\Theta}^1 \cup \bar{\Theta}^1)
\right)\nonumber\\
& \le &
e^{-n/8} + e^{-n/2} + {N\choose\ell} (2e)^{\ell}
e^{-\kappa^2 m \ell/(8 \cdot 6^2)}.
\end{eqnarray}
This ends {\em Step I } of the proof.  To summarize: we have shown that
the exceptional set $\Theta^0$ is of exponentially small measure
provided (\ref{lowerkappa}) holds, and that if, additionally,
(\ref{2cubeball}) is satisfied, then, for $\omega \not\in \Theta^0$,
the quotient space $\tilde{X}_p$  (obtained from $X_p(\omega)$ via the
quotient map $Q$)  contains a well-complemented subspace well-isomorphic to
$W$. To be precise, to arrive at such a conclusion requires optimizing
the estimate (\ref{2singleop-p}) over allowable choices of the parameters
$N, \kappa$; however, we skip it for the moment since an even more subtle
optimization will be performed in {\em Steps II } and {\em III}.

\medskip

{\em Steps II } and {\em III } are very similar as in \cite{ST1},
Proposition 3.1, so we shall outline the main points only, referring
the interested reader to \cite{ST1} for details.

\bigskip
\noindent
{\em Step II.\ The perturbation  argument. \ }
Let $Q$ be an {\em arbitrary } rank $m$ orthogonal projection on
$\Rn$. Denote by $\Theta^Q$
the set given by formally the same formulae as in (\ref{Theta0}) by the
Gaussian operator $\tilde{G} = QG$ for this particular $Q$.
By rotational invariance,  all the properties we
derived for $\Theta^0$  hold also for $\Theta^Q$.
Throughout {\em Step II}, all references to objects defined in {\em Step I}
will implicitly assume that we are dealing with this particular $Q$.

Consider
the exceptional set $ {\Theta}^1 $ defined in (\ref{Theta1}), and
observe that if $\omega \not\in {\Theta}^1$, then
\begin{equation}
   \label{2Omega1inclusion}
   D_{j, p}' \subset D_p \subset 2 B_2^n,
\end{equation}
for every  $j = 1, \ldots, N$.
This is an analogue of (3.26) of \cite{ST1} and the basis for all
the estimates that follow.

\medskip

Let $\omega \not\in \Theta^1 \cup \Theta^Q$ and let
$Q'$ be any  rank $m$ orthogonal projection
such that $\|Q-Q'\| \le \delta$, where $\| \cdot \|$ is the operator
norm with respect to the Euclidean norm $|\cdot |$ and $\delta > 0$
will be specified later. Then,
for some $j$, conditions just slightly weaker than those in
(\ref{Theta'j}) and (\ref{Theta''(j)}) hold with $Q$ replaced by $Q'$.
Namely, there exists $ 1 \le j \le N$ such that, firstly, if $\delta
\le (1/8)\sqrt{m/n}$ then $Q'$ satisfies inclusions analogous to
(\ref{Theta''(j)}) with constants $1/2$ and $2$ replaced by $1/4$ and
$9/4$, respectively (cf. (3.28) of \cite{ST1}); and, secondly, if
$\delta_1 := 4 \delta \sqrt{n/m} \le \kappa/4$ then $Q'$ satisfies
inclusions analogous to (\ref{Theta'j}) with $\kappa$ replaced by
$2\kappa$.  (The former statement is exactly the same as in
\cite{ST1}, and the proof of the latter uses the above inclusion
(\ref{2Omega1inclusion}) instead of (2.26) of \cite{ST1}.)

\smallskip

Finally, set $\delta := 1/(8 \sqrt n)$ (as in \cite{ST1}); then the
condition $\delta \le (1/8)\sqrt{m/n}$ is trivially satisfied, while the
condition $\delta_1 \le \kappa/4$
follows from (\ref{lowerkappa}).  So we can now apply the
previous arguments and conclude {\em Step II}:
if $\omega \not\in \Theta^1 \cup
\Theta^Q$, $\|Q-Q'\| \le \delta $ and
\begin{equation} \label{2cubeball4}
2 \kappa
\leq \frac 1{\sqrt{k}} \cdot \frac14 \sqrt{\frac m{n}} ,
\end{equation}
then the quotient of $X_p$  corresponding to
$Q'$ contains a $2^{1/q}$-complemented
subspace $2^{1/q}$-isomorphic to $W$, namely $Q'E_j$.
We note that (\ref{2cubeball4}) is just slightly stronger than
(\ref{2cubeball}),  and as easy to satisfy.

\medskip

{\em Step III.\ The discretization: a $\delta$-net argument. \ }
Let $\mathcal{Q}$ be a $\delta$-net in the set of rank $ m$ orthogonal
projections on $\R^n$ endowed with the distance given by the operator
norm.  Recall that such a net can be taken with cardinality
$|\mathcal{Q}| \le (C_2/\delta)^{m(n-m)}$, where $C_2$ is a universal
constant (see \cite{ST1}, or directly \cite{S2}).
For our choice of $\delta = 1/(8 \sqrt n)$, this does not exceed
$e^{mn\log{n}}$, at least for sufficiently large $n$.  As in
(\ref{1singleop-p})-(\ref{2singleop-p}), this implies
the measure estimate for our final exceptional set
\begin{eqnarray}
\label{2manyop}
\pp\left(\Theta^1 \cup
\bigcup_{Q\in \mathcal{Q}} \Theta^Q \right) 
& \le & \pp\left(\Theta^1 \cup \bar{\Theta}^1 \cup
\bigcup_{Q\in \mathcal{Q}}
\bigl(\Theta^Q \setminus (\Theta^1 \cup \bar{\Theta}^1)\bigr)\right)\\
&\leq & e^{-n/8} + e^{-n/2} +
e^{mn\log n} \, {{N}\choose{\ell}}
(2e)^\ell
e^{-\kappa^2 m \ell/32} \nonumber
\end{eqnarray}
The first two terms are negligible.
Recall that
$\ell = \lceil N/3 \rceil \ge N/3$, and so the last term in
(\ref{2manyop}) is less than or equal to
$e^{mn\log{n} - \kappa^2 m N/128}$.

\medskip
In conclusion, if  $k$, $\kappa$,  $m$,  $n$ and $N$ satisfy
\begin{equation}
\label{2warunki}
C \sqrt{q} \, (4/\sqrt \pi)
\sqrt{k/m}\, N^{1/q} \le \kappa/12,
\qquad
   256\, mn\log{n} \le  \kappa^2 m N,
\end{equation}
where $C >0$ is the absolute constant related to $C_p$ (see
(\ref{lowerkappa}) and  (\ref{cp_const})),
then the set $\Omega \setminus (\Theta^1 \cup
\bigcup_{Q\in \mathcal{Q}} \Theta^Q)$ has positive measure (in fact,
very close to 1 for large $n$).
If, additionally, (\ref{2cubeball4}) is satisfied, then any $\omega$
from this set induces an $n$ dimensional space $X_p$ whose {\em all }
$m$-dimensional quotients contain a $2^{1/q}$-isomorphic  and
$2^{1/q}$-complemented copy of $W$ (and similarly  with $1+\ep$ in place
of $2^{1/q}$ if  (\ref{2cubeball4}) holds with an additional $\ep$ factor on
the right hand side). Then the assertion of Theorem~\ref{2hmmp} holds for
that particular value of $m$.

\medskip

It remains to ensure that conditions (\ref{2warunki}) and
(\ref{2cubeball4})  are consistent and to
discuss the resulting  restrictions on the dimensions.
It is most convenient to let
$\kappa := (1/8)\sqrt{ m/( n k)}$
so that (\ref{2cubeball4}) holds. Then
the conditions in  (\ref{2warunki}) lead to
\begin{equation}
\label{warunki2p}
k \le c' \min\left\{ \frac{m}{ \sqrt{q n} N^{1/q}},
      {\frac{m N}{n^2 \log n}}\right\},
\end{equation}
where $c' \in (0,1)$ is a universal constant.  Optimizing over $N$
leads to
$$
k \le \frac{c_1 m}{ q^{1/2}\ n^{(4 + q)/(2 + 2q)}
         \ (\log{n})^{1/(1 + q)}}.
$$
which, for $m=m_0$, is just a rephrasing of the hypothesis on $\dim V = \dim
W$ from Theorem \ref{2hmmp}, and holds in the entire range $m_0 \le m \le
n$ if it holds for $m_0$.
It follows that, under our hypothesis, the above construction
can be implemented for each $m$ verifying  $m_0 \le m \le n$.  Moreover,
since the estimates on the probabilities of the exceptional sets
corresponding to different values of $m$ are exponential in $-n$ (as
shown above), the sum of the probabilities involved is small.
Consequently, the construction can be implemented {\em simultaneously}
for all such $m$ with the resulting space satisfying the full
assertion of Theorem \ref{2hmmp} with probability close to $1$.

Finally, we point out that, as it was already alluded to earlier at some
crucial points of the argument,  the $1+\ep$-version of the statement
will follow once our
parameters satisfy (\ref{2warunki}) and the condition analogous to
(\ref{2cubeball4}), with an extra $\ep$ on the right hand side. With the
choice of $\kappa := (\ep/8)\sqrt{ m/( n k)}$, this leads to a version of
(\ref{warunki2p}), which -- after optimizing over $N$ -- gives the same bound
for $k$ as above, but with the constant $c_1$ depending on
$\ep$ rather than being universal.  The rest of the argument is the same.
\qed

\section{ The global saturation}
\label{2global-sec}

\noindent
{\em Proof of Theorem \ref{2hmm-global}{\ }}
Let $W$ be a $k$-dimensional normed space. Identify
$W$ with $\R^k$ in such a way that
$(1/\sqrt k)  B_2^k \subset B_W \subset  B_2^k$.

We use an analogous notation for convex bodies as in the the proof of
Theorem~\ref{2hmmp} (but without the subscript $p$).  In
particular, we set $Z = \ell_1^N (W)$ and we recall that $G =
G(\omega)$ denotes a $n \times Nk$ random matrix with independent
$N(0, 1/n)$-distributed Gaussian entries. We let
$$
K= B_{X(\omega)}  := G (\omega) (B_Z) \subset \R^n.
$$
Recall that for $j=1, \ldots, N$, $F_j$ is the $k$-dimensional
coordinate subspace of $\R^{Nk}$ corresponding to the $j$th
consecutive copy of $W$ in $Z$;  $E_j := {G}(F_j) $,
$K_j := {G}(F_j \cap B_Z)$ and
$K_j' := {G}\left(\spn[F_i: i \ne j] \cap B_Z\right)
= \conv (K_i: i \ne j)$;
next, $D_j := {G}(F_j \cap  B_2^{Nk}) $ and
$
D_j' := \conv  \left(D_j:  i \ne j\right).
$
(The notation $D_j$ has been already used in the proof
of Theorem~\ref{2hmmp}, and the ``$p$-convex''
analogoue of $D_j'$,  namely $D_{j, p}'$,  was
defined in (\ref{djp'}).)

The general structure of the argument is the same as in
Theorem~\ref{2hmmp}: the proof consists of three steps dealing
respectively with analysis of a single rotation, perturbation of a
given rotation and discretization
(for a smoother narrative, here and in what follows we refer
to elements of $O(n)$ -- even those whose determinant is not
$1$ -- as rotations). We will refer extensively to
arguments in Section~\ref{2subspaces} and in \cite{ST1}. As in
Section~\ref{2subspaces}, we shall occasionally assume, as we may, that
$n$ is large.

\bigskip
\noindent {\em Step I.  Probability estimates for a fixed  rotation.\ }
For the time being we fix $u: \R^n \to \R^n$ with $u \in O(n)$.
We shall show that, outside of an exceptional set of $\omega$'s of a
small measure, there is a section of $K + u(K)$ which is $3$-isomorphic to
$B_W$ and $3$-complemented  (or, more precisely,
that the identity on $W$ $3$-factors through the space $(\R^n, K + u(K))$.

\medskip

We shall adopt the following description of the body $K + u(K)$.
Let $B_Z \oplus_\infty B_Z$ be the unit ball of $Z \oplus_\infty
Z$ (i.e.,  $\R^{Nk} \oplus \R^{Nk}$ with the $\ell_\infty$-norm on the
direct sum). Next, consider the Gaussian operator
$G \oplus G: \R^{Nk} \oplus \R^{Nk} \to \R^{n} \oplus \R^{n}$,
acting in the canonical way on the coordinates. Further, define
$[\Id, u]: \R^{n} \oplus \R^{n} \to \R^{n}$ by $[\Id, u] (x_1, x_2) =
x_1 + u x_2$, for $(x_1, x_2) \in \R^{n} \oplus \R^{n}$.  Clearly, we
have $K + u(K) = [\Id, u] (G \oplus G)(B_Z \oplus_\infty B_Z)$.
Instead of $[\Id, u]$ we can equally well use $[u_1, u_2]$, where
$u_1, u_2 \in O(n)$ are two rotations.

The difference between this setup and the scheme of \cite{ST1} is that in
the latter one considers $Q G''(B_Z \oplus_1 B_Z)$, where $G''$ is a $2n
\times 2 Nk$ Gaussian matrix and $Q$ a $\rank n$ orthogonal projection on
$\R^{2n}$.  Both schemes yield quotients of random quotients of $Z \oplus
Z$, with $G
\oplus G$ or $G''$ being the random part and $[u_1, u_2]$ or $Q$ the
nonrandom part.  For the latter one may as well ``rescale" the
dimensions and consider $Q'G(B_Z)$, where $Q'$ is a (nonrandom) rank
$\lfloor n/2 \rfloor$ projection.  The setting in
Section~\ref{2subspaces} is identical, except that we consider
$B_{Z_p} \oplus_p B_{Z_p}$ instead of $B_Z \oplus_1 B_Z$.

\medskip

To define exceptional sets we identify conditions similar to those
in Section~\ref{2subspaces} (or in Section 3 of \cite{ST1}). Recall
that for $E \subset \R^n$, by $P_E$ we denote the orthogonal
projection onto $E$. Now, for $j \in \{1, \ldots, N\}$, and $ 0 <
\kappa <1$ fixed, to be
specified later, we consider the set
\begin{equation}
\label{Omega'j-gl}
\OO_j' :=
\left\{ \omega \in \Omega :
P_{E_j} (D_j' + u(D_j')) \subset
\kappa  B_2^n \right\}.
\end{equation}
These sets are analogous to $\Theta_j'$ in (\ref{Theta'j}), and they will
replace these latter sets in all subsequent definitions.
A similar proof as for
(\ref{2Theta_est}) in Section~\ref{2subspaces},
or (3.23) of \cite{ST1}, shows that
\begin{equation}
\label{both_meas1}
\pp (\OO_j') \ge 1 - \exp (-c_1 \kappa^2 n),
\end{equation}
as long as
\begin{equation}
      \label{k1}
\kappa  \ge C' \sqrt{\max\{k,\log{N}\}/n},
\end{equation}
for appropriate numerical constants $c_1 >0$ and  $C' \ge 1$.
The argument
is again based on Lemma \ref{shrinking}: since $E_j$ is independent of
$D_j'$, we may as well consider it fixed, and then we are exactly in the
setting of the Gaussian part of the Lemma. We just need to majorize
$M^*(B_Z)$ (or, more precisely, just of the unit ball of
$\ell_1^{N-1}(\ell_2^k)$ since the $\ell_2^k$-factor corresponding to $F_j$
does not enter into $D_j'$), which is $O(\sqrt{\max\{k,\log{N}\}/n})$
by reasons similar to -- but simpler than -- those that led to (3.20) of
\cite{ST1} (the calculations sketched in the paragraph containing
(\ref{cp_const})  give a slightly larger majorant, which would also
suffice for our purposes).

\smallskip
Next, for $j=1, \ldots, N$ we  let
\begin{equation}
\label{case2_ball''}
\OO_{j,0}' :=  \left\{\omega \in \Omega :
(1/2) \, (B_2^n \cap E_j) \subset D_j
\right\}.
\end{equation}
Since the condition in (\ref{case2_ball''}) involves only one of the
two inclusions appearing in (\ref{Theta''(j)}), the same argument that led
to (\ref{2bigcube})  (see also (3.16) of \cite{ST1}) gives
\begin{equation}
\label{case2_meas2}
\pp (\OO_{j,0}') \ge  1 - \exp (- n/32).
\end{equation}

\medskip
While in Theorem~\ref{2hmmp} and in
\cite{ST1},  properties analogous to those implicit in the
definitions of the sets $\OO_j'$, $\OO_{j,0}'$
were sufficient to ensure that the quotient
$Q(K)$ contained a well-complemented subspace  well-isomorphic to
$W$, this is not the case in the present context and we need
to introduce additional invariants.

Fix  $\al_0 >0 $ to be specified later (it will be of the order of
$1/k$). Let $\alpha := \tr (\Id - u)/n$, and assume without loss of
generality that $0 \le \alpha \le 1$ (replacing,  if necessary, $u$
by $-u$).
The proof now splits into two cases depending on whether
$\al \ge \al_0$ or $\al < \al_0$.
To clarify the structure of the argument let us mention that, among the sets
$\OO_j'$ and
$\OO_{j,0}'$ defined above, {\em Case} $1^\circ$ will use only the former
ones, while {\em Case} $2^\circ$ will involve both.

\bigskip
\noindent {\em Case} $1^\circ$: Let $\al \ge \al_0$.

\begin{lemma}
\label{case1lemma}
Let $A $ be an $n \times k$ random matrix with independent
$N(0,1/n)$-distributed Gaussian entries.
Let $u \in O(n)$  with ${\tr} u \ge 0$ and
set $\alpha = {\tr}(\Id - u)/n$ ($\in [0,1]$).
Then, with probability
greater than or equal to
$1 - \exp (- c \alpha n + c^{-1} k \log (2/\alpha))$,
the following holds for
all $\ \xi , \zeta  \in \R^k$
\begin{equation}
\label{case1_lowerest}
\left|
A \xi + u A \zeta \right|
\ge c  \alpha ^{1/2} \left(|\xi|^2+|\zeta|^2 \right) ^{1/2}
\ge (c/2) \alpha ^{1/2} \left(|A \xi|^2+|u A \zeta|^2 \right) ^{1/2},
\end{equation}
where $c>0$ is a universal constant.
\end{lemma}

We postpone the proof of the lemma until the end of the section and
continue the main line of the argument.
For $j = 1, \ldots, N$ we let
$$
H_j = E_j + u(E_j).
$$
We shall now use Lemma \ref{case1lemma} for the $n \times k$ matrix
$A= A_j$ formed by the $k$ columns of the matrix $G$ that span  $E_j$.
Denoting by  $\OO_{j,0} $  the subset of
$\Omega$ on which the inequalities (\ref{case1_lowerest}) holds,
we have
\begin{eqnarray}
\label{case1_meas}
\pp (\OO_{j,0}) & \ge & 1 - \exp (- c \alpha n + c^{-1} k \log
(1/\alpha)) \nonumber \\
& \ge & 1 - \exp (- c \alpha_0 n + c^{-1} k \log (1/\alpha_0)).
\end{eqnarray}
Consider the following auxiliary set, closely related to
$\OO_{j,0} $,
\begin{equation}
\label{case1_ball}
\Delta_j  :=  \left\{\omega \in \Omega :
c \alpha ^{1/2} (B_2^n \cap H_j) \subset D_j + u(D_j)
\ \mbox{\rm and}\ \dim H_j = 2k
\right\}.
\end{equation}
An elementary argument shows that the conditions in (\ref{case1_ball})
are equivalent to
``$\left| A \xi + u A \zeta \right|
\ge c  \alpha ^{1/2} \max \{|\xi|, |\zeta| \}$
for all $\ \xi , \zeta  \in \R^k$." Since this is weaker than
the first inequality in (\ref{case1_lowerest}), it follows that
$ \OO_{j,0} \subset \Delta_j$.

\medskip
Our next objective is to show that on $\OO_{j,0}$
\begin{equation}
\label{hjprojection}
|P_{H_j} z| \le  (2/c) \alpha ^{-1/2} \Bigl(|P_{E_j} z|^2 +
|P_{u(E_j)} z|^2 \Bigr)^{1/2},
\end{equation}
for every $z \in \R^n$.

Note that since $E_j$ and $u(E_j)$ are both subspaces of $H_j$, it
is sufficient to assume that $z \in H_j$. Consider the operator $T :
H_j \to E_j \oplus_2 u(E_j)$ given by $T(z) = ( P_{E_j} z, P_{u(E_j)}
z)$ for $z \in H_j$. Then the inequality (\ref{hjprojection}) is
equivalent to  $\|T^{-1}\| \le (2/c) \alpha ^{-1/2}$. On
the other hand, the adjoint operator $T^*: E_j \oplus_2 u(E_j) \to H_j$
is given by $T^*(x, y) = x+y$ for $x \in E_j$ and $y \in u(E_j)$.
Comparing the first and the third terms of (\ref{case1_lowerest})
yields $\|T^{-1}\|=\|{(T^*)}^{-1}\| \le (2/c) \alpha ^{-1/2}$, as required.

\medskip
Finally, consider another  good  set
\begin{equation}
\label{OOmega''j-gl}
{{\OO}}_j'' :=
\left\{ \omega \in \Omega :
P_{u(E_j)}(D_j' + u(D_j')) \subset
\kappa B_2^n \right\}.
\end{equation}
Note that since $u$ is orthogonal,
we clearly have
$P_{u(E_j)} = u P_{E_j} u^*$
(this will be used more than once).  Comparing (\ref{OOmega''j-gl}) with the
definition of $\OO_j'$ (see (\ref{Omega'j-gl})), we deduce from
(\ref{both_meas1}) that
\begin{equation} \label{combined}
\pp({{\OO}}_j'')= \pp({{\OO}}_j') \ge 1 - \exp (-c_1 \kappa^2 n)
\end{equation}

We are now ready to complete the analysis specific to {\em
Case} $1^\circ$.  Let $\omega \in \OO_{j,0} \cap {\OO}_j' \cap {{\OO}}_j''$.
Then, combining  (\ref{hjprojection}) with
the definitions of  $\OO_j'$ and $\OO_j''$
-- i.e., with (\ref{Omega'j-gl})  and (\ref{OOmega''j-gl}) --  we see that,
for all $z \in D_j' + u(D_j')$,
$$
|P_{H_j}z| \le  (2/c) \alpha ^{-1/2} \Bigl(|P_{E_j} z|^2 +
|P_{u(E_j)} z|^2 \Bigr)^{1/2}
\le  (2\sqrt 2/c) \alpha ^{-1/2} \kappa
$$
or, equivalently,
\begin{equation}
\label{OOmega'j-gl}
P_{H_j} (D_j' + u(D_j')) \subset  (2\sqrt 2/c)
    \alpha ^{-1/2}\kappa B_2^n
\end{equation}
As in the previous proofs we will  impose a condition on $\kappa$, namely
\begin{equation}
\label{case1_cond}
(2\sqrt 2/c) \alpha_0 ^{-1/2}\kappa
\le
c  \al_0^{1/2}/\sqrt k.
\end{equation}
Combining this inequality with (\ref{OOmega'j-gl}) and (\ref{case1_ball}),
and recalling that $ \OO_{j,0} \subset \Delta_j$ and that $\al_0 \le \al$,
we are led to
$$
P_{H_j } (D_j' + u(D_j')) \subset 1/\sqrt k \left(D_j + u(D_j)\right).
$$
Finally,  recalling the inclusions between the $K$- and the $D$-sets,
we obtain
$$
P_{H_j } (K_j' + u(K_j'))  \subset K_j + u(K_j.
$$
Consequently, similarly as in the previous proofs (cf.
(\ref{2koniec}), or (3.3) of \cite{ST1}),
$$
P_{H_j}(K + u(K)) \subset \conv \left(K_j + u(K_j),  P_{H_j}(K_j' +
u(K_j'))\right) \subset K_j + u(K_j).
$$
This means that $K_j + u(K_j)$ is a $1$-complemented section of
$K + u(K)$. On the other hand, let us
note that, again by (\ref{case1_ball}), $\dim
H_j = 2k$, which implies that $K_j + u(K_j)$ (thought of as a normed
space) is isometric to $B_W \oplus_\infty B_W$,
thus showing that $H_j \cap (K + u(K))$ is
isometric to $B_W \oplus_\infty B_W$ as well.

We recall that the above conclusion was arrived at under the hypothesis
$\omega \in \OO_{j,0} \cap {\OO}_j' \cap {{\OO}}_j''$. As $j \in \{1, \ldots
, N \}$ was arbitrary, we deduce that under the hypothesis of {\em Case}
$1^\circ$ and the additional assumptions (\ref{k1}) and (\ref
{case1_cond}),  the set $K + u(K)$ admits a
$1$-complemented section isometric to $B_W$ provided that $\omega \in
\bigcup_{j=1}^N (\OO_{j,0} \cap {\OO}_j' \cap {{\OO}}_j'')$.

\bigskip

\noindent
{\em Case} $2^\circ$: Let $\al <  \al_0$.

In this case the operator $u$ is close to the identity operator.
In particular, since $\al = \tr (\Id -u)/n$, we see that the norm
$$
\|\Id - u\|_{HS} = \left( \tr (\Id - u)(\Id - u^*) \right)^{1/2}
= \left(2(n - \tr u) \right)^{1/2} = (2n\al)^{1/2}
$$
is relatively small.
To exploit this property we will need another lemma.

\begin{lemma}
\label{case2lemma}
Let $A $ be an $n \times k$ random matrix with independent
$N(0,1/n)$-distributed Gaussian entries.
Let $T$ be an $n \times n$ matrix,
set $a :=\|T\|_{HS}$ and let $\gamma >0$.  Then,
on a set of
probability larger than or equal to
$1 - \exp \left(- \gamma^2 n/(2 \|T\|^2 )
+ 2 k\right)$,
the following holds for all  $\xi=(\xi_i) \in \R^k$
\begin{equation}
\label{case2_est}
\left| T A \xi \right|
\le 2 \left(\frac{a}{\sqrt n } + \gamma \right) |\xi|.
\end{equation}
\end{lemma}

\noindent
Again, we postpone the proof of the Lemma and continue our argument.
Fix $\gamma >0$,  to be specified later.  For $j = 1, \ldots, N$, let
\begin{equation}
\label{case2_ball}
\OO_{j,0}'' :=  \left\{\omega \in \Omega :
(\Id -u) D_j \subset 2(\sqrt{2\al} + \gamma) B_2^n
\right\}.
\end{equation}

As was the case with Lemma \ref{case1lemma}, we
shall apply the Lemma to the $n \times k$ matrix
$A = A_j$ formed by the $k$ columns of the matrix $G$ that span  $E_j$.
We will also use $T = \Id -u$, so that $\|T\| \le 2$.
Since, in that case, $a/\sqrt n = \sqrt{2\al}$, the inclusion from
(\ref{case2_ball}) is equivalent to the inequality
(\ref{case2_est}) and thus
\begin{equation}
\label{case2_meas}
\pp (\OO_{j,0}'') \ge  1 - \exp (- \gamma^2 n/8 + 2 k).
\end{equation}
The latter expression will be later made very close to 1 by an
appropriate choice of parameters.

\smallskip
Next we shall show that if $j \in \{1, \ldots, N \}$ and
$\omega \in  \OO_j' \cap \OO_{j,0}' \cap \OO_{j,0}''$, then
\begin{equation}
\label{case2_concl}
K_j \subset P_{E_j} (K + u(K)) \subset  3 K_j.
\end{equation}
Clearly, this will imply that the section of $K+ u(K)$ by $E_j$ is
$3$-isomorphic to $K_j$, which in turn is isometric to $B_W$; and
additionally, that it is $3$-complemented.
Consequently,  under the hypothesis of {\em Case} $2^\circ$,
the assertion of {\em Step I } will be shown to hold on the set
$\bigcup_{j=1}^N (\OO_{j,0}' \cap \OO_{j,0}''\cap \OO_j')$.

To show (\ref{case2_concl}), we first point out that if
$B \subset \R^n$ is {\em any } symmetric convex body,  then
$B + u(B) \subset 2B + (\Id-u)(B)$.
We then argue as follows
\begin{eqnarray*}
K + u(K) & \subset & K_j + K_j' + u(K_j) + u(K_j')\\
& \subset &  K_j + D_j' + u(K_j) + u(D_j') \\
& \subset & 2 K_j + (\Id - u)K_j + \left(D_j' + u(D_j') \right)\\
& \subset&  2 K_j +  2(\sqrt{2\al} + \gamma) B_2^n +
\left(D_j' + u(D_j')\right),
\end{eqnarray*}
where the last inclusion is a consequence of (\ref{case2_ball}).
Accordingly
\begin{eqnarray*}
P_{E_j}\left(K + u(K)\right) &\subset&
2 K_j + P_{E_j} \left(D_j' + u(D_j')\right) + 2(\sqrt{2\al} + \gamma)
B_2^n \cap E_j\\
&\subset & 2K_j + (\kappa + 2\sqrt{2\al} + 2\gamma)
(B_2^n \cap E_j),
\end{eqnarray*}
with the last inclusion following from the definition
(\ref{Omega'j-gl}) of  set $\OO_j'$.
By the definition
(\ref{case2_ball''})
of $\OO_{j,0}'$, the second term on the right is contained in
$2 (\kappa + 2\sqrt{2\al} + 2\gamma) D_j$.
Since $\al < \al_0$, it follows that whenever
\begin{equation}
\label{case2_cond}
2 (\kappa + 2\sqrt{2\al_0} + 2\gamma)
\le 1/\sqrt {k} ,
\end{equation}
then
$$
P_{E_j} \left(K + u(K)\right) \subset 2 K_j +
(1/\sqrt {k}) D_j \subset 3 K_j.
$$
We thus obtained the right hand side inclusion in  (\ref{case2_concl});
the left hand side inclusion is trivial. This ends the analysis specific to
{\em Case} $2^\circ$.

\bigskip

Now is the time to choose $\al_0$ and $\gamma$ to satisfy our
restrictions while yielding the optimal concentration in {\em both } cases
under consideration.  The conditions (\ref{k1}), (\ref{case1_cond}) and
(\ref{case2_cond}) can be summarized as
$C' \sqrt{\max\{k,\log{N}\}/n} \le \kappa \le c' \alpha_0/\sqrt k$
and $\max \left\{\sqrt{\al_0}, \gamma\right\}
\le c'/\sqrt k$, for appropriate numerical constants $c' >0$ and $C' \ge 1$.
We choose
$\al_0$, $\gamma$ and $\kappa$ so that
\begin{equation}
\label{both_cond}
\kappa^{1/3} = \sqrt{\al_0} =
\gamma = c'/ \sqrt k
\end{equation}
This choice takes care of all the restrictions except for the lower bound on
$\kappa$, which can be now rephrased as
\begin{equation}
\label{both_cond2}
k \le c \min \{ n^{1/4}, \left(n/\log{N} \right)^{1/3} \},
\end{equation}
for an appropriate numerical constant $c >0$.

\medskip We shall now analyze the estimates on the probabilities
of the good sets contained in (\ref{case2_meas}),
(\ref{case1_meas}) and (\ref{combined}). If $k^2/n$ is sufficiently small,
a condition which is weaker than (\ref{both_cond2}), then the term $2k$ in
the exponent in (\ref{case2_meas}) is of smaller order than the first
term, and so it does not affect the form of the estimate.  The
situation is slightly more complicated in the case of
(\ref{case1_meas}): to absorb the second term in the exponent we
need to know that $k \log{(1/\alpha_0)}$ is sufficiently smaller than
$\alpha_0 n$; , given that $\alpha_0 = O(1/k)$ (cf.
(\ref{both_cond})), this is equivalent to
$$
k \le c'' \ \sqrt{\frac{n}{1+\log{n}}}
$$
for an appropriate numerical constant $c'' >0$. Again, this is a condition
weaker than (\ref{both_cond2}), at least for sufficiently large $n$.  The
probability estimates in question are thus, respectively, of the form
$1 - \exp (- c_3 \gamma^2 n)$, $1 - \exp (- c_2 \alpha_0 n)$
and $1 - \exp (- c_1 \kappa^2 n)$, for appropriate universal constants $c_1,
c_2, c_3 >0$.  Substituting the values for $\al_0$, $\gamma$ and
$\kappa$ defined by (\ref{both_cond})
we get, under the hypothesis (\ref{both_cond2}), the following
minoration
\begin{equation}
\label{both_meas}
\min \{
\pp({\OO}_j'), \pp({{\OO}}_j''), \pp(\OO_{j,0}), \pp(\OO_{j,0}''),
\pp(\OO_{j,0}')
\}
\ge 1 - \exp (- c_0 n/k^3) ,
\end{equation}
again for an appropriate numerical constant $c_0 >0$.  We point out that the
argument above treated just the first four terms under the minimum; for
$\pp (\OO_{j,0}')$  we have the stronger estimate (\ref{case2_meas2}), which
does not require any additional assumptions.

\medskip
We are now ready to conclude {\em Step I}.  Consider the exceptional
set defined by one of two different formulae, depending on whether we are
in {\em Case}
$1^\circ$ or
{\em Case} $2^\circ$. In {\em Case} $1^\circ$ we set
$$
\OO^0 := \Omega \setminus \bigcup_{j=1}^N
(\OO_{j,0} \cap {\OO}_j' \cap {{\OO}}_j'')
$$
(see (\ref{case1_lowerest}) and the paragraph following it,
(\ref{Omega'j-gl}), (\ref{OOmega''j-gl})
for the definitions).
In {\em Case} $2^\circ$ we let
$$
\OO^0 := \Omega \setminus
\bigcup_{j=1}^N
( \OO_{j,0}' \cap \OO_{j,0}''\cap \OO_j'),
$$
(see (\ref{case2_ball''}), (\ref{case2_ball}) and (\ref{Omega'j-gl})
for the definitions).
The argument above shows that for $\omega \not\in \OO^0 $
there is a section of $K + u(K)$ $3$-isomorphic to $B_W$ and
$3$-complemented.

It follows readily from what we have shown up to now that
the sets $\OO^0$ are exponentially small.
For example, by (\ref{both_meas}),
\begin{equation}
\pp\bigl(\Omega \setminus (\OO_{j,0} \cap {\OO}_j' \cap {{\OO}}_j'')\bigr)
\le 3 \exp (- c_0 n/k^3)
\label{singlej}
\end{equation}
for any $j \in \{1, \ldots, N\}$, and identical estimates hold for
exceptional sets relevant to {\em Case} $2^\circ$.
However,  to finalize {\em Step I }  we need to majorize the
probability of $\OO^0 $ much more efficiently.
To this end we argue in the same way as in Section 3 of \cite{ST1}.
We could also follow the argument from Section~\ref{2subspaces} above,
but in the present situation, when we are dealing with the convex
hulls of sets -- such as $K_i$ or $D_j$ -- rather then the $p$-convex
hulls of the same sets, with $p >1$, the latter option would only
add  unnecessary complications. However, for reader's convenience,
we will also include a few comments pertaining to the proof of
Theorem~\ref{2hmmp}.

We first employ the ``decoupling''
procedure based on Lemma 3.2 in \cite{ST1} (which is a special case of
Lemma~\ref{decoupling_lemma} above for a ``0-1'' matrix $A$).
More precisely,
we do need and do have
estimates on conditional probabilities which are obtained in
essentially the same way as there (and are also parallel to the
estimates for $\Theta^0$ earlier in this paper).
Essential use is also made of  the exceptional set
$$
\Omega^1 := \{ \omega : D \not\subset 2 B_2^n \}
$$
(defined in (3.17) of \cite{ST1} and analogous to $\Theta^1$ in
Section~\ref{2subspaces})
and the precise statements
involve $\Omega^1$ and sets related to it.
Again, the key point is that the
linear subspace $E_j$ (resp. $E_j + u(E_j)$) and the sets with which it
is being intersected
(or which are projected onto it)
depend on disjoint blocks of columns of $G$ and hence are independent.
The decoupling procedure  and the estimate from (\ref{both_meas}) lead to
\begin{equation}
\label{small1}
\pp (\OO^0) \le
N e^{-9 n/32} + {{N}\choose{\ell}}
\left(3e^{-c_0 n/k^3}\right)^\ell \le N e^{- 9 n/32} + e^{-c_4 N/k^3},
\end{equation}
where $\ell = \lceil N/3 \rceil $ (cf. (\ref{singlej})). This is almost
identical to (3.25) of
\cite{ST1} (and analogous to (\ref{2singleop-p}) above).
Let us emphasize that the set $\Omega^1$, responsible for the first
term of the estimate, is independent on $u$, and therefore, when
(\ref{small1}) is combined with the $\delta$-net argument in {\em Step
III } below, only the second term will have to be multiplied by the
cardinality of the net.

\bigskip
\noindent {\em Step II.  Stability under small perturbations of the
rotation $u$. \ }
We will now prove that there exists (a not too small) $\delta > 0$
such that if $u \in O(n)$ and $\omega \not\in \OO^0$ (where $\OO^0$ is
defined starting with this particular $u$)
and if  $u' \in O(n)$ with $\|u-u'\| \le \delta$,
then $u'$ and $\omega$ satisfy essentially the same conditions as those
defining $\OO^0$.
As in \cite{ST1} (and analogously as in Section~\ref{2subspaces}
above), this will be shown under an additional assumption, namely that
$\omega \not\in \Omega^1 $
(the definition of $\Omega^1$ was recalled above). It will then follow
that, for any $u'$ as above, the random body $K$ corresponding
to any  $\omega \not\in \Omega^1 \cup \OO^0$ will have the property that
$K + u'(K)$ has a section that is $3$-isomorphic to $B_W$ and
$3$-complemented provided the parameters involved in the construction satisfy
conditions differing from those of {\em Step I } (which, we recall,
were ultimately reduced  just to (\ref{both_cond2})) only by values of the
numerical constants.

We start by pointing out that the condition (\ref{case2_ball''})
does not involve $u$ and so it is trivially stable.
Next, we consider (\ref{Omega'j-gl}) which, while non-trivial,  is easy to
handle. We have
$$
u'(D_j') \subset u(D_j') + 2 \delta B_2^n
$$
(because $\omega \not\in  \Omega^1$) and so if
$\delta \le  \kappa/2$,
we get  (\ref{Omega'j-gl}) for $u'$ in place of $u$,
at the cost of replacing $\kappa$ by $2\kappa$ on the right hand side of the
inequality.

The condition (\ref{case2_ball}) is also  simple:
if  $\omega \not\in  \Omega^1$ and
$\delta \le \gamma $,
and if (\ref{case2_ball}) is satisfied for $u$, then it is clearly
satisfied for $u'$ with the factor $2$ on the right hand side replaced
by $3$.
(Note that this argument works for a general $u$, even though the condition
(\ref{case2_ball}) enters the proof only in {\em Case} $2^\circ$.)

Next we assume that we are in {\em Case} $1^\circ$ and discuss the stability
of $\OO_{j,0}$, defined by inequality (\ref{case1_lowerest}) (where the
matrix $A = A_j$ has been described in the paragraph following
Lemma~\ref{case1lemma}).  We clearly have
\begin{eqnarray*}
\left|
A_j \xi + u' A_j  \zeta \right|
& \ge &\left|
A_j \xi + u A_j  \zeta \right|
- \|u - u'\| \left|A_j \zeta \right| \\
&\ge & c \alpha ^{1/2} \left(|\xi|^2 + |\zeta|^2\right)^{1/2}
- 2 \delta |\zeta|.
\end{eqnarray*}
So if $\delta \le c \alpha ^{1/2}/4$, we get a version of the
first inequality in (\ref{case1_lowerest}) with $u'$ in place of $u$
and $c$ on the right hand side replaced by $c/2$.  The second
inequality follows similarly.  Since (given that we are in {\em Case}
$1^\circ$)  $\al \ge \al_0 $,
we see that the condition on $\delta$ is satistfied when $\delta \le c
\alpha_0 ^{1/2}/4$.


It remains to check the stability of (\ref{OOmega''j-gl}).
Set $R = u' - u$, then $\|R\| \le \delta$ and, using $P_{u'(E_j)}
=u' P_{E_j} {u'}^*$, we obtain
\begin{eqnarray*}
P_{u'(E_j)}(D_j' + u'(D_j')) & = &
  u' P_{E_j}({u'}^*D_j' + D_j')\\
& = & (u +R) P_{E_j}((u^* + R^*) D_j' + D_j')\\
&\subset & (u +R) P_{E_j}\left((u^* D_j'+ D_j') + 2 \delta B_2^n\right)\\
&\subset & u P_{E_j}(u^* D_j' + D_j') + 2  \delta u P_{E_j} B_2^n
+ R P_{E_j}(4 B_2^n)\\
&\subset & u P_{E_j}(u^* D_j' + D_j') + 6  \delta  B_2^n.
\end{eqnarray*}
Since $u P_{E_j} u^* = P_{u(E_j)}$, insisting that $\delta \le \kappa/6$
will guarantee that
$u'$ satifies  the inclusion from (\ref{OOmega''j-gl}) with
$ \kappa$ replaced by $2\kappa$.

Finally, let us remark that the distinction between {\em Cases}
$1^\circ$ and $2^\circ$ is likewise essentially stable under small
perturbations in $u$: the parameter $\al$ is 1-Lipschitz with respect
to the operator norm and so if $\delta$ is less than $1/2$ of the
threshold value $\al_0 = {c'}^2/k$, then the inequalities defining
{\em Case} $1^\circ$ and $2^\circ$ will have to be modified at most by
factor $2$ when passing from $u$ to $u'$ (or vice versa).

Comparing the obtained conditions on $\delta$ we see
that the most restrictive is
$\delta \le \kappa /6 = c''' k^{-3/2}$.
Since, by (\ref{both_cond2}) (and,
ultimately, by the hypothesis of the Theorem), $k$ is at most of the
order of $n^{1/4}$, the appropriate choice of $\delta = O(n^{-3/8})$,
will cover the entire range of possible values of $k$.
This supplies the value of $\delta$ that needs to be used in the
discretization (a $\delta$-net argument) to be
implemented in  {\em Step III } below.

\bigskip
\noindent
{\em Step III.  A discretization argument.\ }
The procedure is fully parallel to that of Section~\ref{2subspaces}: we
introduce a $\delta$-net of $O(n)$, say $\mathcal{U}$, and then
combine  the exceptional sets corresponding to the elements of
$\mathcal{U}$. For the argument to work, it will be sufficient that
the cardinality of $\mathcal{U}$ multiplied by the probability of the
exceptional set corresponding to a particular rotation $u$
(i.e., the second
term at the right hand side of (\ref{small1})) is small.  As is well
known (see, e.g., \cite{S1}, \cite{S2}), $O(n)$ admits, for any
$\delta > 0$, a $\delta$-net (in the operator norm) of cardinality not
exceeding
$(C/\delta)^{\dim {O(n)}}$, where $C$ is a universal constant.  Our
choice of $\delta = O(1/n^\beta)$ (where $\beta = 3/8$,
see the last paragraph of {\em Step II}) leads to the estimate
$$
\log {|\mathcal{U}|} \le O(\beta n^2 (1 + \log {n})) .
$$
For the probability of combined  exceptional sets to be small
it will thus suffice that, for an appropriately chosen $c_5 > 0$,
$$
\beta n^2 (1 + \log {n}) \le c_5 N/k^3
$$
(cf. (\ref{small1})).
Since, as in the argument at the end of {\em Step II}, we may assume
that $k$ is at most of the order of $n^{1/4}$, the condition above may
be satisfied in the entire range of possible values of $k$ with $N =
O(n^{11/4} (1 + \log {n}))$.  Since such a choice implies that ue have
then $\log {N} = O(\log {n})$, the restrictions given by
(\ref{both_cond2})) reduce, at least for large $n$, to $k \le c
n^{1/4}$ -- exactly the hypothesis of the Theorem. \qed

\medskip

To complete the proof of Theorem \ref{2hmm-global}
it remains to prove Lemmas \ref{case1lemma} and \ref{case2lemma}.  The
arguments are fairly straightforward applications of the Gaussian
isoperimetric inequality, or Gaussian concentration, again in the form
given, e.g., in \cite{L}, formula (2.35).

\medskip \noindent {\em Proof of Lemma \ref{case1lemma} }
Fix $\xi, \zeta \in \R^k$ with $|\zeta|^2 + |\xi|^2 =1$ and consider
$ f:= \left|A  \xi + u A  \zeta \right|$ as a function of
the argument $A$.
Then $f$ is $\sqrt{2}$-Lipschitz
with respect to the Hilbert-Schmidt norm.
Therefore Gaussian concentration
inequalities imply that the function $f$
must be strongly concentrated around its expected value ${\E}f$.
Specifically, we get for $t >0$
\begin{equation}
\label{concentration}
\pp (|f - {\E}f| > t) < 2 \exp(-n t^2/4).
\end{equation}
To determine the magnitude of ${\E}f$, we shall first
calculate the second moment.
\begin{eqnarray*}
{\E}f^2 &=& {\E}\left| A \xi + u A \zeta \right|^2 \\
&=& {\E}\left|  A \xi \right|^2
+{\E}\left|  u A \zeta \right|^2
+2 {\E}\langle A  \xi, u A \zeta  \rangle \\
&=& |\xi|^2 + |\zeta|^2
+2 \langle \xi, \zeta  \rangle \ \frac {\tr u} {n},
\end{eqnarray*}
the last equality following, for example, by direct calculation  in
coordinates.
The assumption $|\xi|^2 + |\zeta|^2 = 1$ implies
$|\langle \xi, \zeta  \rangle | \le 1/2$ and so,
recalling  our notation
$\alpha = {\tr(\Id - u)}/{n} = 1- {\tr u}/{n}$,
we deduce that
$$
\alpha \le {\E}f^2 \le 2-\alpha .
$$
Since, by the Khinchine-Kahane inequality, the $L_2$- and
the $L_1$-norm of a Gaussian vector differ at most by factor
$\sqrt{\pi/2}$ (see \cite{LO} for an argment which gives the
optimal value of the constant), it follows that
$$
\ep_1 := \sqrt{2/\pi} \ \sqrt{\alpha}  \le {\E}f
\le \sqrt{2-\alpha} .
$$
Thus choosing $t = \ep_1/3$ in (\ref{concentration}) yields
\begin{equation}
\label{concentration1}
\pp \left(2\ep_1/3 \le |A \xi + u A \zeta| \le
\sqrt{2-\alpha} +\ep_1/3\right)
\ge 1- 2 e^{-n\alpha/(18\pi)}.
\end{equation}
The estimates on $|A \xi + u A \zeta|$ and the associated
probabilities extend  appropriately by homogeneity to any
$\xi, \zeta \in \R^k$.
The next step is now standard: we choose a proper
net in the set  $\{ \xi, \zeta \in \R^k : |\zeta|^2 + |\xi|^2 = 1 \}$
and if the estimates on $|A \xi + u A \zeta|$ hold
simultaneously for all elements of that net, it will follow that
$$
1/3 \sqrt{2/\pi} \ \sqrt{\alpha} \ (|\xi|^2 + |\zeta|^2)^{1/2}
\le |A \xi + u A \zeta|
\le 2 (|\xi|^2 + |\zeta|^2)^{1/2},
$$
for all $\ \xi, \zeta \in \R^k$.  The left hand side inequality
above yields then the first inequality in (\ref{case1_lowerest}).  The
right hand side inequality is a statement formally stronger than the
second inequality in (\ref{case1_lowerest}).

To conclude the argument
we just need to assure the proper resolution
of the net and to check its cardinality. Generally, if a linear
map is bounded from above by $B$ on an $\ep$-net of the sphere,
it is bounded on the entire sphere by $B'=B/(1-\ep)$.
If it is additionally bounded on the net from below by $b$,
then it is bounded from below on the entire sphere by $b'= b-B'\ep$.
If we choose $\ep = \ep_1/6$, then the resulting $B'$ is $< 2$,
and so $b' > 2\ep_1/3 - 2 \ep = \ep_1/3$, as required.
Finally, the  $\ep$-net can be chosen so that its cardinality
is $\le (1 + 2/\ep)^{2k} = (1 + \sqrt{18 \pi/\alpha})^{2k}$, and
so the logarithm of the cardinality is $O(k \log (2/\alpha))$.
Combining this with
(\ref{concentration1}) we obtain an estimate on probability
which is exactly of the type asserted in Lemma  \ref{case1lemma}.
\qed

\bigskip
\noindent {\em Proof of Lemma \ref{case2lemma} }
The argument here is similar to that of  Lemma \ref{case1lemma} but
substantially simpler since we need only an upper estimate.
First, we may assume without loss of generality that $T$ is
diagonal. A direct calculation shows then that
${\E}\left| T A \xi \right|^2 =(\|T\|_{HS}^2/n)\,  |\xi|^2$.
Thus, if $|\xi|=1$, then
${\E}\left| T A \xi \right| \le a/\sqrt{n}$, while the
Lipschitz constant of $| T A \xi |$ (in argument $A$,
with respect to the Hilbert-Schmidt norm) is $\le \|T\|$. It is
now enough to choose a $1/2$-net on the sphere $S^{k-1}$ and argue as
earlier, but paying attention to upper estimates only.  \qed

\bigskip 

\footnotesize
\address

\end{document}